\documentclass[11pt,a4paper]{article}
\usepackage[margin=2.55cm]{geometry}
\usepackage[T1]{fontenc}
\usepackage{lmodern,amsmath,amssymb,amsthm,mathtools,bm,microtype}
\usepackage{booktabs,enumitem,xcolor,array,graphicx,placeins}
\usepackage[colorlinks=true,linkcolor=blue!45!black,citecolor=blue!45!black,urlcolor=blue!45!black]{hyperref}
\hypersetup{pdftitle={Local geometry and limit laws for smooth adapted Wasserstein distances},pdfauthor={Xiaozhen Wang}}
\newcommand{\R}{\mathbb R}
\newcommand{\E}{\mathbb E}
\newcommand{\Pp}{\mathbb P}

\newcommand{\AW}{\mathcal{AW}}
\newcommand{\HH}{\mathcal H}
\newcommand{\QQ}{\mathcal Q}
\newcommand{\DD}{\mathcal D}
\newcommand{\CC}{\mathcal C}

\newcommand{\NN}{\mathcal N}
\newcommand{\supp}{\operatorname{supp}}
\newcommand{\Var}{\operatorname{Var}}
\newcommand{\Cov}{\operatorname{Cov}}
\newcommand{\Law}{\mathcal L}
\newcommand{\BL}{\mathrm{BL}}

\DeclareMathOperator{\diver}{div}
\newtheorem{theorem}{Theorem}[section]
\newtheorem{lemma}[theorem]{Lemma}
\newtheorem{proposition}[theorem]{Proposition}
\newtheorem{corollary}[theorem]{Corollary}
\theoremstyle{definition}
\newtheorem{assumption}[theorem]{Assumption}
\theoremstyle{remark}
\newtheorem{remark}[theorem]{Remark}
\numberwithin{equation}{section}
\setlist{nosep,leftmargin=1.6em}
\allowdisplaybreaks[2]
\title{Limit distributions and bootstrap \\ for smooth adapted Wasserstein distances}
\author{Xiaozhen Wang\thanks{CMAP, \'Ecole polytechnique,
Institut Polytechnique de Paris, Palaiseau, France.}}
\date{November 15, 2026}

\begin{document}
\maketitle

\begin{abstract}
In this paper, we provide limit theorems for the Gaussian-smoothed adapted Wasserstein distances. Let $\widehat\mu_n$ be the empirical measure of independent paths with law $\mu$ on $(\R^d)^T$. Let $P=\mu*\gamma_\sigma$ and $P_n=\widehat\mu_n*\gamma_\sigma$, where $\gamma_\sigma=\NN(0,\sigma^2I_{dT})$ is a Gaussian measure with mean zero and covariance matrix $\sigma^2I_{dT}$. We first assume that $\mu$ is compactly supported. For every $1<p<\infty$ and fixed $d,T$ and $\sigma>0$, we show that
\[
 \sqrt n\,\AW_p(P_n,P)\ \Longrightarrow\ Q_{P,p}(G)^{1/p},
\]
where $G$ is the Gaussian limit of the density fluctuation and $Q_{P,p}$ is the minimum adapted $p$-energy for the linearized continuity equation. We prove full-path bootstrap consistency and the convergence of the first and $p$-th moments. The limit and the bootstrap hold for noncompact laws as well, for all $1<p<\infty$ under appropriate tail and structural conditions. For two-sided bounded density perturbations of strongly log-concave laws, the bandwidth is not restricted for $1<p\le2$. For $p>2$, we give sufficient tail conditions and provide a Gaussian counterexample where root-$n$ tightness fails for small bandwidths. Our proof is based on a two-endpoint local expansion of bicausal transport. As an application, we provide asymptotic error bounds for optimal stopping values under the smoothed model using the bootstrap quantiles.
\end{abstract}

\section{Introduction}\label{sec:introduction}

In multi-period decision problems, comparing two stochastic process laws can be much more involved than just measuring their proximity in the ordinary Wasserstein metric on the path space $(\R^d)^T$. This is because the decision made at time $t$ can only use information available up to that time. Hence, two path laws can be very close with respect to the ordinary Wasserstein metric $W_p$ but can still differ substantially in their conditional distributions of future states given the past. This in turn can change the values of optimal stopping or hedging problems. To compare multi-stage models along their filtrations, Pflug \cite{Pflug2009} and subsequently Pflug and Pichler \cite{PP2012} developed the concepts of nested distributions and nested distance. More recently, Backhoff-Veraguas, Bartl, Beiglb\"ock, and Eder established the equivalence of various adapted weak topologies \cite{BBBEtop} and used adapted Wasserstein distances $\AW_p$ (defined through bicausal couplings) to prove Lipschitz stability results in mathematical finance \cite{BBBE}.

When we turn from population stability to estimating a distribution from data, the sensitivity to conditional laws becomes a problem. For example, suppose we observe $n$ independent trajectories $A^{(1)},\dots,A^{(n)}$ from a law $\mu$ with an atomless first marginal. Then with probability one, the starting points $A_1^{(1)},\dots,A_1^{(n)}$ are all distinct. Under the raw empirical measure $\widehat\mu_n$, observing the first coordinate therefore identifies the sampled path and reveals its entire future. As a consequence, $\widehat\mu_n$ need not converge to $\mu$ in $\AW_p$ as $n\to\infty$. Backhoff, Bartl, Beiglb\"ock, and Wiesel \cite{BBBW} addressed this difficulty by constructing adapted empirical measures through partitions of the state space and deriving quantitative convergence rates.

A natural alternative is regularization by convolution smoothing: one convolves both measures with a Gaussian kernel $\gamma_\sigma=\NN(0,\sigma^2 I_{dT})$ before computing the transport distance. For Wasserstein distances, this smoothing approach was introduced by Goldfeld and Greenewald \cite{GG2020} and developed for general $p$ in \cite{NGK2021}, while the corresponding central limit theorems and bootstrap theory were established in \cite{GKNR}. In the multi-period setting, having a smooth joint density is only half the story. One also has to control how the conditional densities vary with the past trajectory. Recent work has addressed several facets of this issue: the bandwidth-dependent behavior of smoothed empirical measures in $\AW_p$ \cite{Hou}, quantitative bounds comparing adapted and classical transport metrics \cite{BLPW}, and $\AW_p$ projections restoring martingale constraints in empirical pricing \cite{BWZZ}. Recently, Larsson, Park, and Wiesel \cite[Theorem~4]{LPW} showed that at fixed $\sigma>0$ and under a bandwidth-dependent quadratic exponential-moment condition, $\E[\AW_p(\widehat\mu_n*\gamma_\sigma,\mu*\gamma_\sigma)]$ decays at rate $n^{-1/2}$. Whether the rescaled error admits a limiting distribution was left open in \cite[Section~1.3]{LPW}.

We address this question. Fix $\sigma>0$, $1<p<\infty$, and integers $d,T\ge 1$, and put $P=\mu*\gamma_\sigma$ and $P_n=\widehat\mu_n*\gamma_\sigma$ with densities $\rho$ and $\rho_n$ on $\R^{dT}$. We study the asymptotic behavior of $\sqrt{n}\,\AW_p(P_n,P)$ and prove, under the hypotheses on $\mu$ given below, that
\begin{equation*}
    \sqrt{n}\,\AW_p(P_n,P)\ \Longrightarrow\ N_{P,p}(G)
    \qquad\text{as }n\to\infty.
\end{equation*}
In this limit, $G$ is the centered Gaussian field arising from the empirical density fluctuation $\sqrt{n}(\rho_n-\rho)$, and $N_{P,p}=Q_{P,p}^{1/p}$ is an adapted tangent norm defined on an appropriate space of zero-mass signed densities. Here, $Q_{P,p}(h)$ is the least $L^p(P)$ kinetic energy of a velocity field $v=(v_1,\ldots,v_T)$ solving the linearized continuity equation $-\sum_{t=1}^T\nabla_t\cdot(\rho v_t)=h$ while respecting causality, in the sense that each $v_t$ is a function of $x_{1:t}$ alone. In the $L^p(P)$--$L^{p'}(P)$ dual variational problem, with $p'=p/(p-1)$, this restriction replaces the test gradients $\nabla_t\phi(X)$ from static smooth Wasserstein limits \cite{GKNR} by their conditional expectations $\E_P[\nabla_t\phi(X)\mid\mathcal F_t]$ \cite{ABGHP,AKPR}.

Since both $N_{P,p}$ and the covariance of $G$ depend on the unknown law $\mu$, direct evaluation of the limit distribution requires estimating these quantities. We prove that the ordinary bootstrap consistently estimates the law of $N_{P,p}(G)$ once we resample entire trajectories $A^{(i)}$ instead of drawing each time coordinate separately; full-path resampling retains the dependence across dates inside the resampled path law. In Section~\ref{sec:application}, we use these bootstrap quantiles together with the Lipschitz stability of optimal stopping to obtain asymptotic error bounds for stopping values under the smoothed model $P$, separately from the deterministic smoothing bias between $P$ and $\mu$. Because the bootstrap distance $\AW_p(P_n^*,P_n)$ involves two random perturbations of $P$ at the same time, the proof relies on a deterministic two-endpoint expansion 
\[
 |\varepsilon|^{-1}\AW_p(P+\varepsilon h,P+\varepsilon j)
 \longrightarrow N_{P,p}(j-h)
 \qquad\text{as }\varepsilon\to0,
\]
valid under our assumptions whenever $h$ and $j$ are Gaussian smoothings of bounded zero-mean tilts of $\mu$. In static optimal transport, the upper bound along a linear density interpolation $\rho_s$ is usually proved by rescaling a velocity field by $\rho/\rho_s$, which in the adapted setting can introduce dependence on future coordinates and spoil adaptedness. Instead, we transport $P+\varepsilon h$ along smooth block-triangular ODE flows and bound the leftover error to $P+\varepsilon j$ through conditional $L^{p'}$-Poincar\'e inequalities.

When $\mu$ has compact support, the limit law, full-path bootstrap consistency, and convergence of the first and $p$-th moments hold for all $1<p<\infty$ and all $\sigma>0$. At $p=2$, the optimal velocity is linear in $h$, which makes $Q_{P,2}(G)$ a weighted series of independent chi-square random variables; for $p\neq 2$, the tangent energy is nonquadratic, as can be seen in a two-period example in Appendix~\ref{app:fourpoint}. Without compact support, the cases $1<p\le 2$ and $p>2$ behave quite differently. For $1<p\le 2$, the limit law and bootstrap hold when $d=1$ and $\mu$ satisfies a quadratic exponential-moment condition relative to $\sigma$, and in any dimension $d\ge 1$ when $\mu=f\mu_0$ is a two-sided bounded density perturbation of a strongly log-concave law $\mu_0(da)\propto e^{-V_0(a)}da$ with $\nabla^2V_0\succeq\kappa I_{dT}$; for such perturbations every fixed $\sigma>0$ works, and the first and $p$-th moments converge too. For $p>2$, both noncompact classes still admit the same limit and bootstrap once stronger bandwidth-dependent exponential-moment conditions are imposed, whereas log-concavity alone does not suffice at small bandwidths: already for $d=T=1$ and $\mu=\NN(0,s^2)$, $\sqrt{n}\,W_p(P_n,P)$ diverges to $+\infty$ in probability whenever $\sigma^2<(p-2)s^2$, because the factor $\rho^{1-p}$ in the tangent energy grows too fast in the tails. Sections~\ref{sec:mainresults} and~\ref{sec:noncompact} state the compact and noncompact theorems, while Section~\ref{sec:application} develops the optimal stopping bound alongside a scalar energy formula and a numerical experiment. The proofs are given in Sections~\ref{sec:expansion} and~\ref{sec:noncompact-estimates}, and Appendix~\ref{app:examples} contains the spectral formulas, further examples, and numerical details.

\section{Main results}

\subsection{Compactly supported path laws}\label{sec:mainresults}

\subsubsection{Setting and adapted tangent norm}

Fix $1<p<\infty$, a state dimension $d\ge 1$, a time horizon $T\ge 1$, and a smoothing scale $\sigma>0$, and put $N=dT$. These parameters remain fixed as $n\to\infty$. A trajectory in $(\R^d)^T$ is written $x=(x_1,\ldots,x_T)$, with prefix $x_{1:t}=(x_1,\ldots,x_t)$ for $1\le t\le T$ and $x_{1:0}=\varnothing$. We write $|\cdot|$ for the Euclidean norm and set $|x|_{\ell^p}=(\sum_{t=1}^T|x_t|^p)^{1/p}$.

Let $\mu$ be a Borel probability measure on $(\R^d)^T$, and let $A^{(1)},\ldots,A^{(n)}$ be independent copies of a random path $A=(A_1,\ldots,A_T)\sim\mu$. Write $\gamma_\sigma=\NN(0,\sigma^2 I_N)$, with the same notation for its density and lower-dimensional marginals. From the empirical measure $\widehat\mu_n=n^{-1}\sum_{i=1}^n\delta_{A^{(i)}}$, we form the smoothed laws
\[
 P=\mu*\gamma_\sigma,\qquad P_n=\widehat\mu_n*\gamma_\sigma,
 \qquad \rho=\frac{dP}{dx},\qquad \rho_n=\frac{dP_n}{dx}.
\]
Writing $e_a(x)=\gamma_\sigma(x-a)=(2\pi\sigma^2)^{-N/2}e^{-|x-a|^2/(2\sigma^2)}$ for $a,x\in\R^N$, we have $\rho(x)=\E_\mu[e_A(x)]$ and $\rho_n(x)=n^{-1}\sum_{i=1}^n e_{A^{(i)}}(x)$. For expectations involving both $A$ and $X$, we use the coupling
\[
 X=A+\sigma Z\sim P,\qquad Z\sim\NN(0,I_N),\quad Z\perp A,
\]
with the same convention for other latent laws indicated by a subscript.

Let $\mathcal F_t=\sigma(X_1,\ldots,X_t)$ be the canonical filtration, with $\mathcal F_0$ trivial.  For two path laws $P$ and $Q$ on $(\R^d)^T$, let $\Pi(P,Q)$ be the set of their couplings, with coordinate processes $(X,Y)$. The set $\Pi_{\rm bc}(P,Q)$ of \emph{bicausal couplings} consists of all $\pi\in\Pi(P,Q)$ such that, for every $1\le t\le T$,
\[
\begin{aligned}
 \Law_\pi(Y_{1:t}\mid X_{1:T})
 &=\Law_\pi(Y_{1:t}\mid X_{1:t}),\\
 \Law_\pi(X_{1:t}\mid Y_{1:T})
 &=\Law_\pi(X_{1:t}\mid Y_{1:t}),
\end{aligned}
\qquad \pi\text{-a.s.}
\]
The adapted $p$-Wasserstein distance is defined by
\begin{equation}
 \AW_p(P,Q)^p=
 \inf_{\pi\in\Pi_{\rm bc}(P,Q)}
       \int\sum_{t=1}^T|x_t-y_t|^p\,d\pi(x,y).
 \label{eq:cost}
\end{equation}
For unsmoothed laws $\nu,\widetilde\nu$, we also put $\AW_p^{(\sigma)}(\nu,\widetilde\nu):=\AW_p(\nu*\gamma_\sigma,\widetilde\nu*\gamma_\sigma)$. When $T=1$, $\AW_p$ reduces to the ordinary $p$-Wasserstein distance $W_p$ on $\R^d$.

\begin{assumption}[Compact support]\label{ass:compact}
The support of $\mu$ is contained in $[-R,R]^N$ for some $R<\infty$.
\end{assumption}

For $1\le t\le T$, let $P_{1:t}$ be the marginal law of $X_{1:t}$ under $P$, with density
$\rho_{1:t}(x_{1:t})=\int\rho(x_{1:t},x_{t+1:T})\,dx_{t+1:T}$. Set $\rho_{1:0}=1$ and let $k_t^P(x_{1:t-1},x_t)=\rho_{1:t}(x_{1:t})/\rho_{1:t-1}(x_{1:t-1})$ be the one-step conditional density. We omit $P$ when clear and also write $k_t^P(z)$ for the measure $k_t^P(z,y)\,dy$. Adapted $L^p$ velocity fields $v=(v_1,\dots,v_T)$, with $v_t=v_t(x_{1:t})$, form the space
\[
 \mathcal A_{P,p}=\bigoplus_{t=1}^T
       L^p(\rho_{1:t}(x_{1:t})\,dx_{1:t};\R^d),
 \qquad
 \|v\|_{P,p}^p=\sum_{t=1}^T\int_{\R^N}
                   |v_t(x_{1:t})|^p\rho(x)\,dx.
\]
Write $\nabla_t$ and $\nabla_t\cdot$ for the gradient and divergence in the $t$-th block $x_t\in\R^d$. For $h\in L^1(\R^N)$ with $\int h(x)\,dx=0$, define
\begin{equation}
 \begin{split}
 Q_{P,p}(h)&=
 \inf_{\substack{v\in\mathcal A_{P,p}\\
                 -\sum_{t=1}^T\nabla_t\cdot(\rho v_t)=h}}
       \sum_{t=1}^T\int_{\R^N}|v_t(x_{1:t})|^p\rho(x)\,dx,\\
 N_{P,p}(h)&=Q_{P,p}(h)^{1/p}.
 \end{split}
 \label{eq:Q}
\end{equation}
The constraint in \eqref{eq:Q} is understood distributionally: $\int\phi h=\sum_{t=1}^T\int\nabla_t\phi\cdot v_t\rho$ for every $\phi\in C_c^\infty(\R^N)$. We set $\inf\varnothing=+\infty$. For $Q_{P,p}(h)<\infty$, the feasible set is closed and affine, so reflexivity and strict convexity of $\mathcal A_{P,p}$ give a unique minimizer. On finite-energy perturbations, $N_{P,p}$ is the norm induced by the quotient of $\mathcal A_{P,p}$ by the closed divergence kernel $\{w\in\mathcal A_{P,p}:\sum_{t=1}^T\nabla_t\cdot(\rho w_t)=0\}$.

The functional \eqref{eq:Q} also admits a dual representation. Let $p'=p/(p-1)$ and identify the dual of $\mathcal A_{P,p}$ with $\mathcal A_{P,p'}=\bigoplus_{t=1}^T L^{p'}(P_{1:t};\R^d)$. Since $v_t(X_{1:t})$ is $\mathcal F_t$-measurable, every $\phi\in C_c^\infty(\R^N)$ and $v\in\mathcal A_{P,p}$ satisfy
\[
 \sum_{t=1}^T\int_{\R^N}\nabla_t\phi\cdot v_t\rho
 =\sum_{t=1}^T\E_P\bigl[\E_P[\nabla_t\phi(X)\mid\mathcal F_t]\cdot v_t(X_{1:t})\bigr].
\]
The divergence kernel is therefore the common null space of these pairings. By reflexivity, its annihilator in $\mathcal A_{P,p'}$ is the norm closure of the conditional gradients $(\E_P[\nabla_t\phi\mid\mathcal F_t])_{t=1}^T$, $\phi\in C_c^\infty(\R^N)$. Putting
\[
 \mathcal E_{P,p'}(\phi)=
 \sum_{t=1}^T\E_P\left|\E_P[\nabla_t\phi\mid\mathcal F_t]\right|^{p'},
 \qquad \phi\in C_c^\infty(\R^N),
\]
quotient duality and Young's inequality yield
\begin{align*}
 N_{P,p}(h)&=
 \sup_{\substack{\phi\in C_c^\infty(\R^N)\\
                     \mathcal E_{P,p'}(\phi)\le1}}\int\phi h,
 \\
 Q_{P,p}(h)&=
 \sup_{\phi\in C_c^\infty(\R^N)}
       \left\{p\int\phi h-(p-1)\mathcal E_{P,p'}(\phi)\right\}.
\end{align*}

\subsubsection{Limit distribution and bootstrap}

Under the hypotheses of the limit theorems below, Section~\ref{sec:densityspaces} specifies a separable Hilbert space $\HH$ of signed densities on which $h\mapsto\int h$ is continuous, $\E_\mu\|e_A\|_\HH^2<\infty$, and $N_{P,p}$ is continuous on $\HH_0=\{h\in\HH:\int h=0\}$. Let $G$ be a centered $\HH_0$-valued Gaussian element with covariance operator
\[
 C_Gu=\E_\mu\big[\langle e_A-\rho,u\rangle_\HH(e_A-\rho)\big],
 \qquad u\in\HH.
\]
The law of $N_{P,p}(G)$ depends only on the covariance kernel $\Cov_\mu(e_A(x),e_A(y))$ and \eqref{eq:Q}, not on the choice of $\HH$.

Conditionally on $\mathcal D_n=\sigma(A^{(1)},\ldots,A^{(n)})$, let $I_1,\ldots,I_n$ be i.i.d.\ uniform on $\{1,\ldots,n\}$, and form the full-path bootstrap laws
\[
 \widehat\mu_n^*=n^{-1}\sum_{j=1}^n\delta_{A^{(I_j)}},\qquad
 P_n^*=\widehat\mu_n^*\ast\gamma_\sigma,\qquad
 \rho_n^*=\frac{dP_n^*}{dx}.
\]
Write $\Law^*,\Pp^*,\E^*,\Var^*$ for conditional law, probability, expectation, and variance given $\mathcal D_n$. On a metric space $(E,d_E)$, set
\[
 d_{\BL,E}(\nu,\widetilde\nu)=
 \sup_{\|f\|_\infty\le1,\,\operatorname{Lip}_{d_E}(f)\le1}
 \left|\int f\,d\nu-\int f\,d\widetilde\nu\right|,
\]
with $d_{\BL}:=d_{\BL,\R}$; on $\R$, the Kolmogorov distance is $\sup_{x\in\R}|\nu((-\infty,x])-\widetilde\nu((-\infty,x])|$.

\begin{theorem}[Compact-support limit and bootstrap]\label{thm:compact}
Under Assumption~\ref{ass:compact}, we have $Q_{P,p}(G)<\infty$ almost surely, and as $n\to\infty$,
\begin{equation}
\begin{aligned}
 \sqrt n\,\AW_p(P_n,P)&\Longrightarrow N_{P,p}(G),
\\
 n^{p/2}\AW_p(P_n,P)^p&\Longrightarrow Q_{P,p}(G).
\end{aligned}
\label{eq:mainlimit}
\end{equation}
The full-path bootstrap satisfies
\begin{equation}
 d_{\BL}\!\left(
 \Law^*\bigl(\sqrt n\,\AW_p(P_n^*,P_n)\bigr),\,\Law\bigl(N_{P,p}(G)\bigr)
 \right)\longrightarrow0
 \quad\text{in probability}.
 \label{eq:bootstrap}
\end{equation}
If $\mu$ is not a point mass, $N_{P,p}(G)$ has a continuous distribution function and \eqref{eq:bootstrap} also holds in Kolmogorov distance; otherwise both $\AW_p(P_n,P)$ and $\AW_p(P_n^*,P_n)$ vanish identically. Moreover, for $r\in\{1,p\}$,
\begin{equation}
 \E\bigl[\bigl(\sqrt n\,\AW_p(P_n,P)\bigr)^r\bigr]
       \longrightarrow\E\bigl[N_{P,p}(G)^r\bigr].
 \label{eq:compactmoments}
\end{equation}
\end{theorem}

\begin{remark}[The quadratic limiting law]\label{rem:quadratic}
For $p=2$, the minimum-energy map $\mathcal V_P:\HH_0\to\mathcal A_{P,2}$ in \eqref{eq:Q} is bounded and linear. Let $\lambda_1\ge\lambda_2\ge\cdots\ge0$ be the eigenvalues of the trace-class operator $\mathcal V_PC_G\mathcal V_P^*$ on $\mathcal A_{P,2}$ (counted with multiplicity, padded with zeros if finite-rank). For i.i.d.\ $Z_j\sim\NN(0,1)$,
\[
 n\,\AW_2(P_n,P)^2
 \Longrightarrow Q_{P,2}(G)
 \overset{d}{=}\sum_{j=1}^{\infty}\lambda_j Z_j^2,
\]
with $\sum_{j=1}^\infty\lambda_j<\infty$; see Appendix~\ref{sec:spectrum} for the proof and trace identity. For $p\neq 2$, the tangent energy $Q_{P,p}$ is nonquadratic.
\end{remark}

\subsection{Laws with unbounded support}\label{sec:noncompact}

When $\mu$ has unbounded support on $(\R^d)^T$, we consider two classes: scalar path laws ($d=1$) with a quadratic exponential moment, and vector-valued laws ($d,T\ge 1$) whose density with respect to a strongly log-concave reference law $\mu_0$ is bounded above and away from zero.

\begin{assumption}[Scalar exponential moment]\label{ass:scalar}
For $d=1$, we suppose that $\E_\mu[e^{\beta|A|^2}]<\infty$ for some $\beta>0$.
\end{assumption}

\begin{assumption}[Strongly log-concave reference law]\label{ass:structured}
There exist $\kappa>0$ and a potential $V_0\in C^2(\R^N;\R)$ such that $\nabla^2V_0(a)\succeq\kappa I_N$ for every $a\in\R^N$. Define the reference probability law by $\mu_0(da)=\mathcal Z_0^{-1}e^{-V_0(a)}\,da$, where $\mathcal Z_0=\int_{\R^N}e^{-V_0(a)}\,da\in(0,\infty)$. We assume that $\mu$ is absolutely continuous with respect to $\mu_0$ and that its density $f=d\mu/d\mu_0$ satisfies
\begin{equation}
 0<c_f\le f\le C_f<\infty
 \qquad \mu_0\text{-almost everywhere},
 \label{eq:referenceclass}
\end{equation}
for some constants $c_f,C_f$.
\end{assumption}

For integers $m\ge 1$, define $\theta_m$ and, when $p>2$, $\Theta_{p,m}$ by
\begin{equation}
\begin{gathered}
 \theta_m=\frac{2m-1+\sqrt{4m^2+12m-7}}4,
 \qquad m\le\theta_m<m+\tfrac12,\\
 \Theta_{p,m}=
 \frac{p(m+1)-2+\sqrt{[p(m+1)-2]^2+8pm}}4
 <\frac{p(m+1)}2.
\end{gathered}
\label{eq:tailthresholds}
\end{equation}
These thresholds arise from the weighted estimates in Sections~\ref{sec:optimizedtails} and~\ref{gpstat:sec-high-p}.

\begin{theorem}[Scalar laws with an exponential moment]\label{thm:scalar}
Suppose Assumption~\ref{ass:scalar} holds with $\beta\sigma^2>\theta_T$ if $1<p\le 2$, or $\beta\sigma^2>\Theta_{p,T}$ if $p>2$. Then \eqref{eq:mainlimit}--\eqref{eq:bootstrap} hold (with Remark~\ref{rem:quadratic} at $p=2$). If $\mu$ is not a point mass, $N_{P,p}(G)$ has a continuous distribution function and \eqref{eq:bootstrap} also holds in Kolmogorov distance; otherwise both $\AW_p(P_n,P)$ and $\AW_p(P_n^*,P_n)$ vanish identically.
\end{theorem}

By the upper bounds in \eqref{eq:tailthresholds}, the assumptions of Theorem~\ref{thm:scalar} hold whenever
\[
 \E_\mu\exp\!\left(\frac{T+\tfrac12}{\sigma^2}|A|^2\right)<\infty
 \quad(1<p\le 2),
 \qquad
 \E_\mu\exp\!\left(\frac{p(T+1)}{2\sigma^2}|A|^2\right)<\infty
 \quad(p>2).
\]

\begin{theorem}[Laws with a strongly log-concave reference]\label{thm:reference}
Suppose Assumption~\ref{ass:structured} holds. For $1<p\le2$, let $\sigma>0$ be arbitrary. For $p>2$, fix $\sigma>0$ and assume that $\E_\mu[e^{\beta|A|^2}]<\infty$ for some $\beta>0$ satisfying $\beta\sigma^2>\Theta_{p,dT}$. Then, as $n\to\infty$,
\[
 \sqrt n\,\AW_p(P_n,P)\Longrightarrow N_{P,p}(G).
\]
The limit has a continuous distribution function, and the conditional full-path bootstrap law $\Law^*(\sqrt n\,\AW_p(P_n^*,P_n))$ converges to $\Law(N_{P,p}(G))$ in Kolmogorov distance, in probability. For $p=2$, the limit is described in Remark~\ref{rem:quadratic}. Moreover, if $1<p\le2$, then
\[
 \E\!\left[\bigl(\sqrt n\,\AW_p(P_n,P)\bigr)^r\right]
 \longrightarrow
 \E\!\left[N_{P,p}(G)^r\right],
 \qquad r\in\{1,p\}.
\]
\end{theorem}

The following Gaussian example shows that, for $p>2$, strong log-concavity alone does not ensure root-$n$ tightness at every bandwidth.

\begin{proposition}[Gaussian obstruction]\label{prop:failure}
Let $d=T=1$, $p>2$, and $\mu=\NN(0,s^2)$ with $s>0$. For every fixed bandwidth satisfying $0<\sigma^2<(p-2)s^2$, we have
\begin{equation}
 \sqrt n\,W_p(\widehat\mu_n*\gamma_\sigma,\mu*\gamma_\sigma)
       \longrightarrow+\infty\quad\text{in probability}.
 \label{eq:failure}
\end{equation}
In particular, the sequence $\{\sqrt n\,W_p(P_n,P)\}_{n\ge1}$ is not tight.
\end{proposition}

Proposition~\ref{prop:failure} is proved in Section~\ref{sec:statistics}.

\subsection{Two-endpoint local expansion}

Our proof rests on a local expansion of $\AW_p$ along bounded latent tilts:
\begin{equation}
 \DD=\{h_g=(g\mu)*\gamma_\sigma:
                    g\in L^\infty(\mu),\ \E_\mu g=0\},
 \qquad P_\varepsilon^h=(\rho+\varepsilon h)\,dx.
 \label{eq:tilts}
\end{equation}
Since $|h_g|\le\|g\|_\infty\rho$, $P_\varepsilon^{h_g}$ is a probability law for $|\varepsilon|\|g\|_\infty\le 1$.

\begin{proposition}[Two-endpoint local expansion]\label{prop:twopoint}
Suppose that one of the following conditions holds:
\begin{enumerate}[label=(\roman*)]
\item Assumption~\ref{ass:compact}. \item Assumption~\ref{ass:scalar}, with
$\beta\sigma^2>\theta_T$ if $1<p\le2$, or
$\beta\sigma^2>\Theta_{p,T}$ if $p>2$.
\item Assumption~\ref{ass:structured}, with the additional
requirement, when $p>2$, that
$\E_\mu[e^{\beta|A|^2}]<\infty$ for some $\beta>0$
satisfying $\beta\sigma^2>\Theta_{p,dT}$.
\end{enumerate}
Here $\theta_m$ and $\Theta_{p,m}$ are defined
in \eqref{eq:tailthresholds}. Then
\begin{equation}
 \lim_{\varepsilon\to0}
 \frac{\AW_p(P_\varepsilon^h,P_\varepsilon^j)}{|\varepsilon|}
 =N_{P,p}(j-h),\qquad h,j\in\DD.
 \label{eq:twopoint}
\end{equation}
\end{proposition}

The two-endpoint form handles the random pair $(P_n^*,P_n)$ in the bootstrap; Section~\ref{sec:statistics} extends \eqref{eq:twopoint} from $\DD$ to empirical perturbations via Hilbert-norm comparisons.

\section{Application and numerical experiment}\label{sec:examples}\label{sec:application}

The scalar formula below makes the limiting energy directly computable. We then apply the bootstrap to an optimal stopping problem and test the resulting calibration in a two-date financial model. The spectral calculation and further examples are collected in Appendix~\ref{app:examples}.

\subsection{Computable tangent energies}\label{sec:computable}

In dimension $d=1$, the continuity equation determines the adapted velocity uniquely whenever a solution in $\mathcal A_{P,p}$ exists, so the infimum in \eqref{eq:Q} disappears. To write out the resulting energy for $T=2$, label the coordinates on $\R^2$ by $(x,z)$, and put $\rho_1(x)=\rho_{1:1}(x)$, $k_x(z)=\rho(x,z)/\rho_1(x)$, and $K(x,z)=\int_{-\infty}^z k_x(s)\,ds$. Given a zero-mass perturbation $h$, define
\[
 h_1(x)=\int_\R h(x,z)\,dz,\qquad
 H_h(x)=\int_{-\infty}^x h_1(s)\,ds,\qquad
 A_h(x,z)=\int_{-\infty}^z h(x,s)\,ds.
\]
Let $v=(v_1,v_2)\in\mathcal A_{P,p}$ satisfy $-\partial_x(\rho v_1)-\partial_z(\rho v_2)=h$. Since both fluxes are integrable, integrating over $z\in\R$ gives $-\partial_x(\rho_1 v_1)=h_1$ and hence $v_1(x)=-H_h(x)/\rho_1(x)$. Substituting $\rho v_1=-H_h k_x$ back into the equation and integrating in $z$ over $(-\infty,z]$ yields $\rho v_2=-A_h+Kh_1+H_h\partial_x K$. The tangent energy is therefore
\begin{equation}
 Q_{P,p}(h)=
 \int_\R\frac{|H_h(x)|^p}{\rho_1(x)^{p-1}}\,dx
 +\iint_{\R^2}\frac{|A_h-Kh_1-H_h\partial_xK|^p}{\rho^{p-1}}\,dx\,dz.
 \label{eq:two-scalar}
\end{equation}
Here the term $-H_h\partial_x K=\rho_1 v_1\partial_x K$ arises from the variation of the conditional CDF $K(x,\cdot)$ in $x$. Repeating this elimination for $T\ge 3$ scalar periods leads to the recursive velocity \eqref{eq:scalarrecursion}. Once $d\ge 2$, the adapted divergence no longer pins down $v$: passing to $dT$ scalar coordinates still produces a feasible field in \eqref{eq:Q}, though not in general the minimizer of the block energy. Under multiparameter conditional laws, the compatibility of local quadratic-energy minimizing velocities is studied in \cite{WWZ2026}.

\subsection{Optimal stopping and bootstrap calibration}

Let $Q$ be a path law on $(\R^d)^T$ equipped with the canonical filtration $\mathcal F_t$, and let $g_t:\R^d\to\R$ ($1\le t\le T$) be measurable discounted payoffs with a common Lipschitz constant $L<\infty$ and finite expectation under $Q$. Writing $\mathcal T(Q)$ for the set of $\{1,\dots,T\}$-valued $\mathcal F_t$-stopping times, we consider the optimal stopping value
\[
 \mathsf V(Q)=\sup_{\tau\in\mathcal T(Q)}\E_Q[g_\tau(X_\tau)].
\]
For a Bermudan put on the first coordinate $x^{(1)}$ with strike $K>0$ and per-period rate $r\ge 0$, the discounted payoff
\[
 g_t(x)=e^{-rt}(K-e^{x^{(1)}})_+,
 \qquad u_+=\max\{u,0\},
\]
is $L$-Lipschitz with $L=K$.

\begin{proposition}[Stability of stopping values]
\label{prop:stopping-stability}
For any path laws $P,Q$ on $(\R^d)^T$ with finite $p$-th moments and $1\le p<\infty$,
\[
 |\mathsf V(P)-\mathsf V(Q)|\le L\,\AW_p(P,Q),
\]
under the cost convention in \eqref{eq:cost}.
\end{proposition}
\begin{proof}
Let $\pi\in\Pi_{\rm bc}(P,Q)$ be a bicausal coupling of $(X,Y)$, and let $\tau\in\mathcal T(P)$. Because $\pi$ is causal from $Y$ to $X$, the conditional law of $\tau$ given $Y$ defines a randomized stopping rule for $Y$. Its expected payoff under $Q$ is a mixture of payoffs of ordinary stopping times and is therefore at most $\mathsf V(Q)$. By H\"older's inequality,
\[
 \E_P[g_\tau(X_\tau)]-\mathsf V(Q)
 \le \E_\pi\bigl[g_\tau(X_\tau)-g_\tau(Y_\tau)\bigr]
 \le L\Bigl(\E_\pi\sum_{t=1}^T|X_t-Y_t|^p\Bigr)^{1/p}.
\]
Taking the supremum over $\tau\in\mathcal T(P)$ and then the infimum over $\pi\in\Pi_{\rm bc}(P,Q)$ gives $\mathsf V(P)-\mathsf V(Q)\le L\,\AW_p(P,Q)$. Swapping $P$ and $Q$ yields the reverse bound.
\end{proof}

Fix a significance level $\alpha\in(0,1)$, and define the full-path bootstrap statistic and its conditional $(1-\alpha)$-quantile by
\[
 T_n^*=\sqrt n\,\AW_p(P_n^*,P_n),\qquad
 q^*_{n,1-\alpha}=\inf\{z\in\R:\Pp^*(T_n^*\le z)\ge 1-\alpha\}.
\]
Let $F_\infty(z)=\Pp(N_{P,p}(G)\le z)$. Its $(1-\alpha)$-quantile $q_{1-\alpha}=\inf\{z:F_\infty(z)\ge 1-\alpha\}$ is unique when $F_\infty(q_{1-\alpha}-\varepsilon)<1-\alpha<F_\infty(q_{1-\alpha}+\varepsilon)$ for every $\varepsilon>0$.

\begin{corollary}[Bootstrap error radius for stopping values]
\label{cor:stopping-bootstrap}
Suppose the assumptions of Theorem~\ref{thm:compact}, Theorem~\ref{thm:scalar}, or Theorem~\ref{thm:reference} hold, and suppose that $N_{P,p}(G)$ is nondegenerate with a unique $(1-\alpha)$-quantile. Then
\[
 \liminf_{n\to\infty}\Pp\left(
 |\mathsf V(P_n)-\mathsf V(P)|
 \le \frac{L\,q^*_{n,1-\alpha}}{\sqrt n}
 \right)\ge 1-\alpha.
\]
\end{corollary}
\begin{proof}
Kolmogorov convergence of $\Law^*(T_n^*)$ to $F_\infty$ together with uniqueness of $q_{1-\alpha}$ yields $q^*_{n,1-\alpha}\to q_{1-\alpha}$ in probability, so
\[
 \Pp\bigl(\sqrt n\,\AW_p(P_n,P)\le q^*_{n,1-\alpha}\bigr)\longrightarrow 1-\alpha.
\]
By Proposition~\ref{prop:stopping-stability}, the event $\sqrt n\,\AW_p(P_n,P)\le q^*_{n,1-\alpha}$ implies $|\mathsf V(P_n)-\mathsf V(P)|\le n^{-1/2}L\,q^*_{n,1-\alpha}$, which gives the result.
\end{proof}

In computation, $q^*_{n,1-\alpha}$ is replaced by the empirical quantile $\widehat q^*_{n,1-\alpha}$ from $B_n$ conditionally independent bootstrap draws. Given the data, the empirical bootstrap CDF has pointwise variance at most $1/(4B_n)$. When $B_n\to\infty$, a finite-grid argument combined with continuity of $F_\infty$ and \eqref{eq:bootstrap} yields uniform CDF convergence in joint probability. Uniqueness of $q_{1-\alpha}$ then implies $\widehat q^*_{n,1-\alpha}\to q_{1-\alpha}$ in probability, so Corollary~\ref{cor:stopping-bootstrap} holds with $\widehat q^*_{n,1-\alpha}$ in place of $q^*_{n,1-\alpha}$ (see Appendix~\ref{app:numerics} for the resampling convention used in our experiment).

\paragraph{The fixed-bandwidth target.}
Write $P_\sigma=\mu*\gamma_\sigma$ and $P_{n,\sigma}=\widehat\mu_n*\gamma_\sigma$. Corollary~\ref{cor:stopping-bootstrap} calibrates the stopping value at a specified smoothing scale $\sigma>0$. If the observed state is $Y=A+\sigma Z$, decisions use the history of $Y$, and exercise pays $g_t(Y_t)$, then $\mathsf V(P_\sigma)$ is directly the value of that observation model; when $\sigma$ is instead chosen as a modeling tolerance, $\mathsf V(P_\sigma)$ serves as a regularized benchmark whose sampling error is calibrated by the bootstrap.

For inference on the unsmoothed value $\mathsf V(\mu)$, let $b_\sigma=|\mathsf V(P_\sigma)-\mathsf V(\mu)|$. Any upper bound $\bar b_\sigma\ge b_\sigma$ adds directly to the statistical radius:
\[
 |\mathsf V(P_{n,\sigma})-\mathsf V(\mu)|
 \le |\mathsf V(P_{n,\sigma})-\mathsf V(P_\sigma)|+\bar b_\sigma.
\]
Under uniformly $W_1$-Lipschitz one-step conditional kernels for $\mu$, the bandwidth bound in \cite[Theorem~5.1]{Hou} together with Proposition~\ref{prop:stopping-stability} at $p=1$ gives $b_\sigma\le L\,\AW_1(P_\sigma,\mu)=O(\sigma)$. On the other hand, for $p>2$ there is a direct obstruction to shrinking the bandwidth $\sigma_n\downarrow 0$: in the Gaussian setting of Proposition~\ref{prop:failure}, fix any $\sigma_0$ satisfying $0<\sigma_0^2<(p-2)s^2$. Once $\sigma_n\le\sigma_0$, contraction of $W_p$ under common Gaussian convolution gives
\[
 W_p(\widehat\mu_n*\gamma_{\sigma_0},\mu*\gamma_{\sigma_0})
 \le W_p(\widehat\mu_n*\gamma_{\sigma_n},\mu*\gamma_{\sigma_n}),
\]
so $\sqrt n\,W_p(\widehat\mu_n*\gamma_{\sigma_n},\mu*\gamma_{\sigma_n})\to+\infty$ in probability along every sequence $\sigma_n\downarrow 0$.

\subsection{A semi-synthetic Bermudan experiment}\label{sec:experiment}

We illustrate the limit law and the full-path bootstrap in a two-date Cox--Ross--Rubinstein binomial model \cite{CRR1979}. Taking the Cboe VIX close of $14.95$ on 31 December 2025 \cite{CBOEVIX} as the annual volatility $\nu=0.1495$, we set $S_0=100$, strike $K=105$, half-year steps $\Delta=1/2$, and an annual interest rate of $3\%$. With
\[
 a=\nu\sqrt\Delta,\qquad U=e^a,\qquad D=e^{-a},\qquad
 \pi_+=\frac{e^{0.03\Delta}-D}{U-D},
\]
each period moves up with risk-neutral probability $\pi_+$, making the discounted stock price a martingale under the latent law $\mu$. In log-relative coordinates, the four paths and their probabilities are
\begin{equation}
\begin{split}
 (a_i)_{i=1}^4&=\bigl((-a,-2a),(-a,0),(a,0),(a,2a)\bigr),\\
 (w_i)_{i=1}^4&=\bigl((1-\pi_+)^2,(1-\pi_+)\pi_+,
                         \pi_+(1-\pi_+),\pi_+^2\bigr).
\end{split}
\label{eq:experiment-law}
\end{equation}
We set $\sigma=0.65a$ and take the discounted Bermudan put payoffs $g_t(x)=e^{-0.03t\Delta}(105-100e^x)_+$ ($t=1,2$), which share the Lipschitz constant $L=105$. Because $\mu$ has finite support, Theorem~\ref{thm:compact} applies for both $p=2$ and $p=4$.

For any probability weights on the four paths in \eqref{eq:experiment-law}, Gaussian smoothing yields second-period conditionals that are stochastically increasing in the first coordinate. The common-uniform Rosenblatt coupling is therefore optimal, allowing us to compute $\AW_p$ directly by one-dimensional quantile quadrature; Appendix~\ref{app:numerics} gives the proof, the quadrature details, and the stopping-value integration.

\begin{figure}[tbp]
\centering
\includegraphics[width=0.98\linewidth]{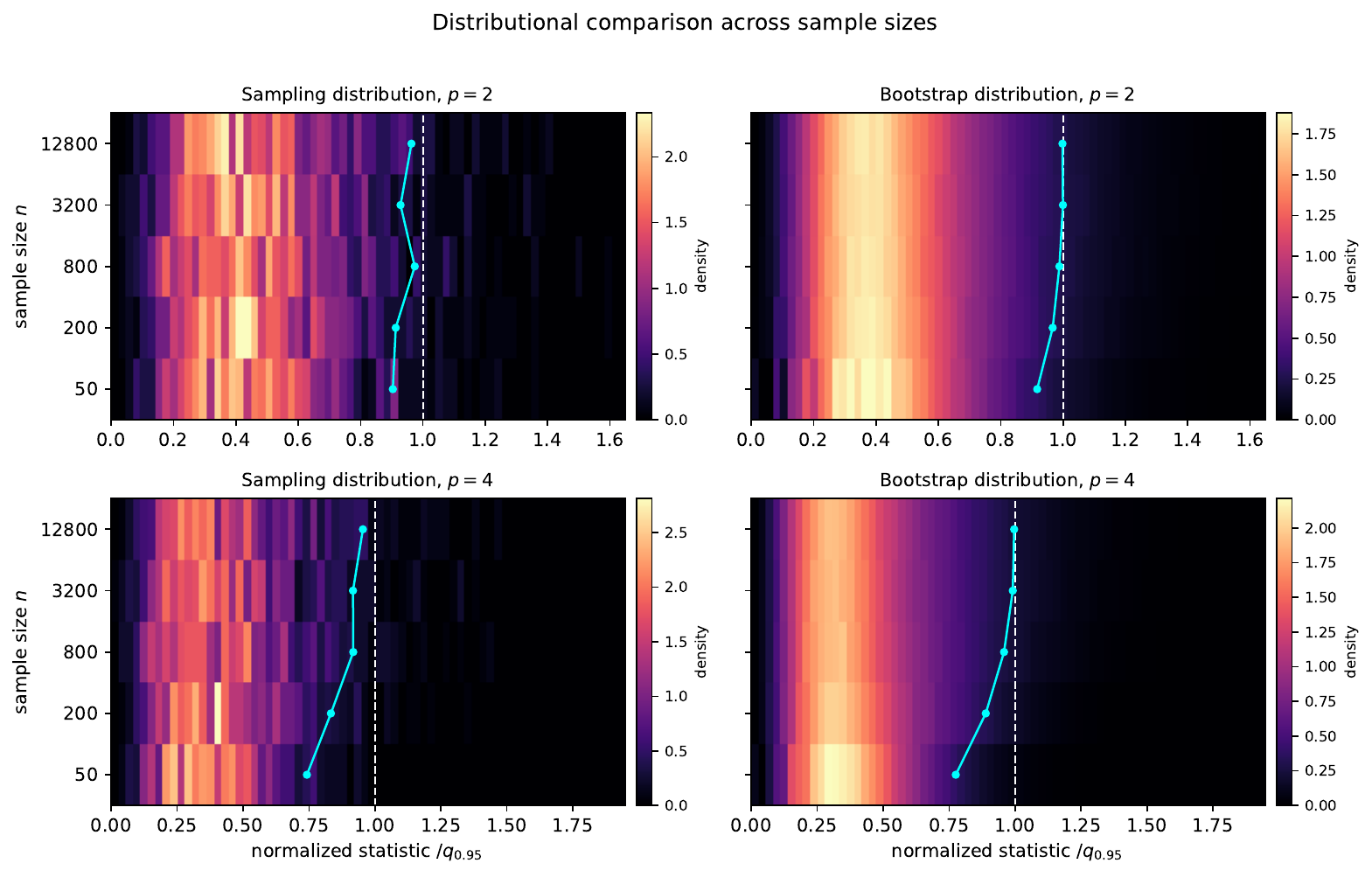}
\caption{Sampling and pooled full-path bootstrap distributions in the model \eqref{eq:experiment-law} for $n\in\{50,200,800,3200,12800\}$. Horizontal axes are normalized by the limiting $0.95$-quantile $q_{0.95}$ from \eqref{eq:two-scalar} (dashed vertical line); the cyan polyline traces the empirical $0.95$-quantile at each $n$.}
\label{fig:distheat}
\end{figure}

\begin{figure}[tbp]
\centering
\includegraphics[width=0.98\linewidth]{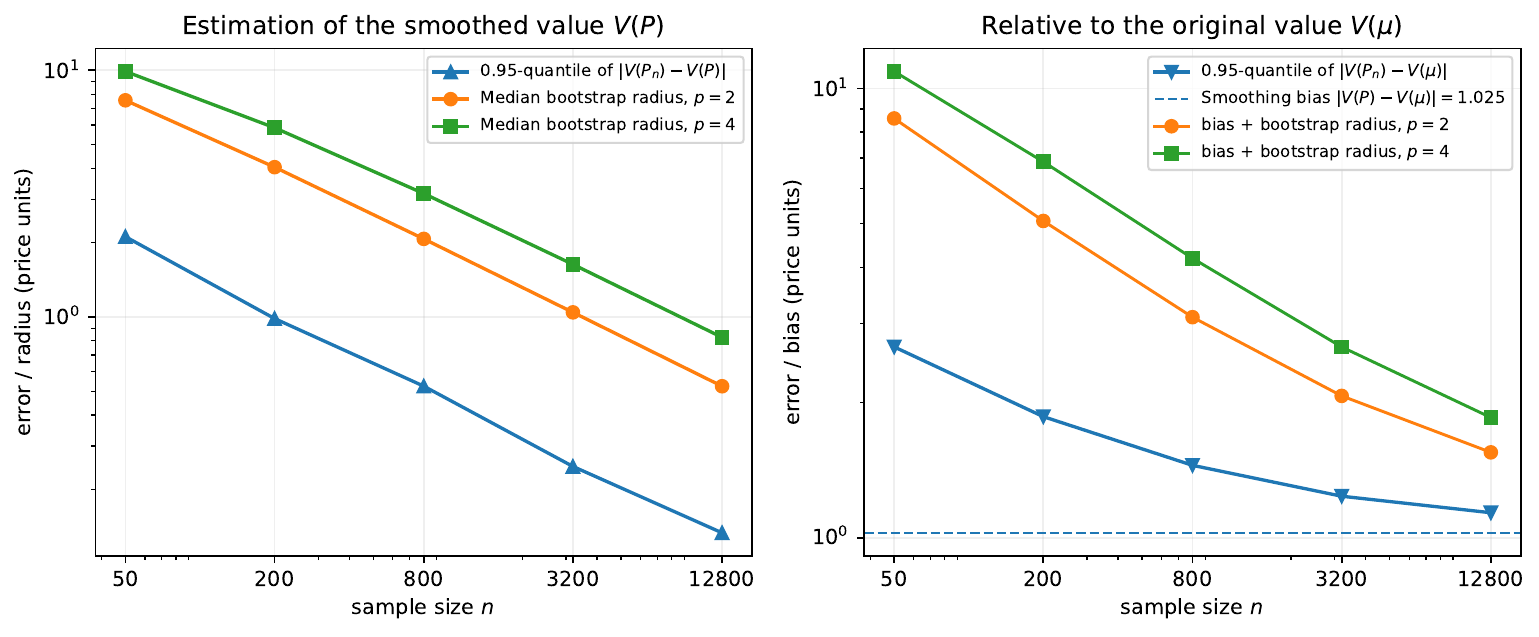}
\caption{Bermudan valuation errors and bootstrap radii (in price units). Left: errors around the smoothed value $\mathsf V(P)$, showing the empirical $0.95$-quantile of $|\mathsf V(P_n)-\mathsf V(P)|$ (blue) and the median bootstrap radii $L\widehat q^*_{n,0.95}/\sqrt n$ for $p=2$ (orange) and $p=4$ (green). Right: the same experiment relative to $\mathsf V(\mu)$, with the empirical $0.95$-quantile of $|\mathsf V(P_n)-\mathsf V(\mu)|$ (blue), the fixed smoothing bias $|\mathsf V(P)-\mathsf V(\mu)|$ (dashed), and this bias added to the corresponding bootstrap radii (orange and green).}
\label{fig:valuation}
\end{figure}

For each $n\in\{50,200,800,3200,12800\}$, we simulate $600$ independent datasets from \eqref{eq:experiment-law} and draw $249$ full-path bootstrap resamples per dataset, recomputing $\AW_p(P_n^*,P_n)$ for each resample. Figure~\ref{fig:distheat} compares the sampling distributions of $\sqrt n\,\AW_p(P_n,P)$ with the pooled bootstrap distributions of $\sqrt n\,\AW_p(P_n^*,P_n)$, normalized by the limiting $0.95$-quantile $q_{0.95}$ computed from \eqref{eq:two-scalar}. Figure~\ref{fig:valuation} compares the resulting stopping-value errors and bootstrap radii around $\mathsf V(P)$ and $\mathsf V(\mu)$.

From \eqref{eq:two-scalar}, the limiting $0.95$-quantiles are $q_{0.95}\approx 0.5666$ for $p=2$ and $q_{0.95}\approx 0.8947$ for $p=4$. Even at $n=50$, where the finite-sample upper tail sits below its asymptotic limit (most visibly for $p=4$), the bootstrap tracks the sampling distribution closely. As $n$ increases, the normalized $0.95$-quantiles of both distributions converge to $1$, and across all five sample sizes the empirical coverage of the nominal $95\%$ distance bound lies between $94.2\%$ and $96.3\%$ (Table~\ref{tab:experiment-coverage}).

Backward induction gives $\mathsf V(\mu)=7.94165$, while numerical integration under the smoothed model gives $\mathsf V(P)=8.96675$, yielding a fixed smoothing bias of $1.02511$. At $n=12800$, the median bootstrap value radii are $0.5235$ ($p=2$) and $0.8261$ ($p=4$), compared with an empirical $0.95$-quantile of $0.1332$ for $|\mathsf V(P_n)-\mathsf V(P)|$. Across all sample sizes, the bootstrap intervals cover $\mathsf V(P)$ in all $600$ trials; adding the smoothing bias $|\mathsf V(P)-\mathsf V(\mu)|$ gives the corresponding error bound around $\mathsf V(\mu)$, as shown in the right panel of Figure~\ref{fig:valuation}.

\FloatBarrier

\section{Proofs of the main results}\label{gpgeo:section}\label{sec:expansion}

Throughout this section, we write $q=p'$ and $L^q(f)=L^q(f(x)\,dx)$ for a positive density $f$. A map or vector field on $(\R^d)^T$ is \emph{block triangular} if its $t$-th component depends only on $x_{1:t}$. For a smooth map $F$, $DF$ denotes its Jacobian matrix and $\|DF\|$ its Euclidean operator norm.

\subsection{Transport comparisons and adapted divergence}

For a path law $P$ with positive density $\rho$, our comparison bounds use the uniform conditional $L^q$-Poincar\'e inequality
\begin{equation}
 \|\phi-\textstyle\int\phi k_t(z,\cdot)\|_{L^q(k_t(z,\cdot))}
 \le C_q\|\nabla\phi\|_{L^q(k_t(z,\cdot))},
 \qquad z\in\R^{d(t-1)}.
 \label{gpgeo:conditionalpoincare}
\end{equation}

\begin{lemma}[Ordinary transport comparison]
\label{gpgeo:ordinarycomparison}
Let $f,g$ be probability densities on $\R^d$ with finite $p$-th moments. If $f$ satisfies the $L^q$-Poincar\'e inequality with constant $C_q$, then
\[
 W_p(f\,dx,g\,dx)
 \le pC_q\left\|\frac{g-f}{f}\right\|_{L^p(f\,dx)}.
\]
\end{lemma}
\begin{proof}
For $\|(g-f)/f\|_{L^p(f)}<\infty$, H\"older's inequality and the $L^q(f)$-Poincar\'e inequality give, for every $\phi\in C_c^\infty(\R^d)$,
\[
 \left|\int\phi(g-f)\right|
 \le C_q\left\|\frac{g-f}{f}\right\|_{L^p(f)}
                 \|\nabla\phi\|_{L^q(f)}.
\]
Hahn--Banach extension and $L^q(f;\R^d)$ duality produce $v\in L^p(f;\R^d)$ solving $-\diver(fv)=g-f$ with $\|v\|_{L^p(f)}\le C_q\|(g-f)/f\|_{L^p(f)}$. Along the nonlinear interpolation
\[
 f_s=(1-\vartheta(s))f+\vartheta(s)g,
 \qquad \vartheta(s)=1-(1-s)^p,\quad s\in[0,1],
\]
the flux $J_s=\vartheta'(s)fv$ satisfies $\partial_s f_s+\diver J_s=0$ with $f_s\ge(1-\vartheta(s))f$ for $s<1$. Since
\[
 \int_0^1\frac{\vartheta'(s)^p}{(1-\vartheta(s))^{p-1}}\,ds
                 \int|v|^pf
 =p^p\int|v|^pf,
\]
the continuity-equation bound for $W_p$ \cite[Chapter~8]{AGS} yields the result.
\end{proof}

\begin{lemma}[One-sided adapted comparison]
\label{gpgeo:onesided}
Let $P,Q$ be path laws with densities $\rho,\eta$, prefix densities $\rho_{1:t},\eta_{1:t}$, and finite $p$-th moments. Suppose $P$ satisfies \eqref{gpgeo:conditionalpoincare},
\[
 W_p(k_t(z,\cdot),k_t(z',\cdot))
       \le L|z-z'|_{\ell^p},
\]
and $\eta\ge c\rho$ for some $c>0$. Then
\[
 \AW_p(P,Q)
 \le C(p,T,C_q,L,c)
        \left\|\frac{\eta-\rho}{\rho}\right\|_{L^p(P)}.
\]
\end{lemma}
\begin{proof}
Set $u=\eta/\rho$ and $u_t=\E_P[u\mid X_{1:t}]$ ($u_0=1$). At any shared history $z=x_{1:t-1}$, we have $k_t^Q(z,\cdot)/k_t^P(z,\cdot)=u_t/u_{t-1}$. Since $u_{t-1}\ge c$, applying Lemma~\ref{gpgeo:ordinarycomparison} at history $z$ and using conditional Jensen's inequality gives
\begin{align*}
 \int \eta_{1:t-1}(z)
       W_p(k_t^P(z,\cdot),k_t^Q(z,\cdot))^p\,dz
 &\le(pC_q)^p\E_P[u_{t-1}^{1-p}|u_t-u_{t-1}|^p]\\
 &\le 2^p(pC_q)^pc^{1-p}\|u-1\|_{L^p(P)}^p.
\end{align*}
Coupling $(X,Y)$ inductively across $t=1,\dots,T$ via measurable optimal (or $\varepsilon$-optimal) one-step couplings of $k_t^P(X_{1:t-1},\cdot)$ and $k_t^Q(Y_{1:t-1},\cdot)$, we insert $k_t^P(Y_{1:t-1},\cdot)$ into the $W_p$ triangle inequality to see that $S_t=\sum_{i\le t}\E|X_i-Y_i|^p$ obeys
\[
 S_t\le(1+2^{p-1}L^p)S_{t-1}
      +2^{2p-1}(pC_q)^pc^{1-p}\|u-1\|_{L^p(P)}^p.
\]
Iterating from $t=1$ to $T$ proves the lemma.
\end{proof}

When $|\nabla_i\log k_t|\le B$ for all $i<t$, we have $|\log(k_t(z',y)/k_t(z,y))|\le B\sum_{i<t}|z_i-z_i'|$, so applying Lemma~\ref{gpgeo:ordinarycomparison} across $m\to\infty$ equal subdivisions of the segment $[z,z']$ yields
\begin{equation}
 W_p(k_t(z,\cdot),k_t(z',\cdot))
 \le pC_qB\sum_{i<t}|z_i-z_i'|
 \le pC_qBT^{1/q}|z-z'|_{\ell^p}.
 \label{gpgeo:historylipschitz}
\end{equation}

On $\mathcal A_{P,p}$, define the distributional adapted divergence operator by
\[
 D_\rho v=-\rho^{-1}\sum_{t=1}^T\nabla_t\cdot(\rho v_t),
\]
with graph domain $\{v\in\mathcal A_{P,p}:D_\rho v\in L^p(P)\}$ equipped with the graph norm $\|v\|_{P,p}+\|D_\rho v\|_{L^p(P)}$. In this notation, the constraint in \eqref{eq:Q} reads $D_\rho v=h/\rho$.

\begin{lemma}[Solvability of adapted divergence]
\label{gpgeo:rightinverse}
Suppose $\rho$ and its prefix densities $\rho_{1:t}$ are positive and smooth, \eqref{gpgeo:conditionalpoincare} holds, and the history scores $\nabla_i\log k_t$ ($i<t$) are bounded. Then for every $s\in L^p(P)$ with $\E_P s=0$, there is a unique minimal-norm velocity field $v\in\mathcal A_{P,p}$ solving $D_\rho v=s$, and
\[
 \|v\|_{P,p}\le C\|s\|_{L^p(P)},
 \qquad
 N_{P,p}(h)\le C\|h/\rho\|_{L^p(P)}\quad\Bigl(\int h=0\Bigr).
\]
\end{lemma}
\begin{proof}
Put $s_t=\E_P[s\mid X_{1:t}]$ with $s_0=0$. Given $v_1,\dots,v_{t-1}$, consider the $t$-th prefix equation
\[
 -\rho_{1:t}^{-1}\nabla_t\cdot(\rho_{1:t}v_t)
 =s_t-s_{t-1}+\sum_{i<t}v_i\cdot\nabla_i\log k_t=:b_t.
\]
Differentiating $\int k_t(z,y)\,dy=1$ in $z_i$ ($i<t$) shows that $\E_P[\nabla_i\log k_t\mid X_{1:t-1}]=0$, so $\E_P[b_t\mid X_{1:t-1}]=0$. By \eqref{gpgeo:conditionalpoincare} and H\"older's inequality, every $\phi\in C_c^\infty(\R^{dt})$ satisfies
\[
 \left|\int b_t\phi\,dP_{1:t}\right|
 \le C_q\|b_t\|_{L^p(P_{1:t})}\|\nabla_t\phi\|_{L^q(P_{1:t})}.
\]
Hahn--Banach extension and $L^q(P_{1:t};\R^d)$ duality yield $v_t\in L^p(P_{1:t};\R^d)$ solving the prefix equation with $\|v_t\|_{L^p(P_{1:t})}\le C_q\|b_t\|_{L^p(P_{1:t})}$. Induction over $t=1,\dots,T$ produces a solution $v\in\mathcal A_{P,p}$ to $D_\rho v=s_T=s$ satisfying the norm bound; reflexivity and strict convexity of $\mathcal A_{P,p}$ then give a unique element of minimal norm in the closed affine space $D_\rho^{-1}(\{s\})$.
\end{proof}

\begin{lemma}[Adapted partial-divergence core]
\label{gpgeo:core}
If, in addition to the assumptions of Lemma~\ref{gpgeo:rightinverse}, the future scores $\nabla_t\log(\rho/\rho_{1:t})$ are bounded for $1\le t\le T$, then every $v\in\mathcal A_{P,p}$ with $D_\rho v=s\in L^p(P)$ can be approximated by smooth fields $v^m=(v_1^m,\dots,v_T^m)$ with $v_t^m\in C_c^\infty(\R^{dt};\R^d)$ such that
\[
 \|v^m-v\|_{P,p}\longrightarrow0,
 \qquad
 \|D_\rho v^m-s\|_{L^p(P)}\longrightarrow0.
\]
\end{lemma}
\begin{proof}
Since $\rho v_i\in L^1(\R^N)$, testing $D_\rho v=s$ against $\psi(x_{1:t})\chi_M(x_{t+1:T})$ for a smooth future cutoff $\chi_M$ with $|\nabla\chi_M|\le C/M$ and letting $M\to\infty$ shows that the partial divergences $D_tv_t:=-\rho_{1:t}^{-1}\nabla_t\cdot(\rho_{1:t}v_t)$ satisfy
\[
 D_tv_t
 =s_t-s_{t-1}+\sum_{i<t}v_i\cdot\nabla_i\log k_t \in L^p(P_{1:t}),
 \qquad s_t=\E_P[s\mid X_{1:t}],\quad s_0=0.
\]
Multiplying $v_t$ by a smooth prefix cutoff $\chi_M(x_{1:t})$ with $|\nabla_t\chi_M|\le C/M$ gives $\chi_M v_t\to v_t$ and $D_t(\chi_Mv_t)=\chi_MD_tv_t-\nabla_t\chi_M\cdot v_t\to D_tv_t$ in $L^p(P_{1:t})$ as $M\to\infty$. At fixed $M$, both $F=\rho_{1:t}\chi_Mv_t$ and $\nabla_t\cdot F$ lie in $L^p(\R^{dt})$ with compact support; convolving $F$ with a standard mollifier $\zeta_\delta$ on $\R^{dt}$ and dividing by $\rho_{1:t}$ gives fields $v_t^{M,\delta}=\rho_{1:t}^{-1}(F*\zeta_\delta)\in C_c^\infty(\R^{dt};\R^d)$ converging to $\chi_M v_t$ in $L^p(P_{1:t})$ and in $D_t$-norm as $\delta\downarrow 0$. Because the future scores $\nabla_t\log(\rho/\rho_{1:t})$ are bounded in the decomposition
\[
 D_\rho v=\sum_{t=1}^T\left(
 D_tv_t-v_t\cdot\nabla_t\log\frac\rho{\rho_{1:t}}\right),
\]
taking a diagonal sequence in $(M,\delta)$ across $t=1,\dots,T$ yields the claim.
\end{proof}

\begin{remark}[Gaussian multipliers]
\label{gpgeo:gaussianmultiplier}
The hypotheses of Lemmas~\ref{gpgeo:rightinverse} and~\ref{gpgeo:core} hold for any density of the form $\widetilde\rho=m\gamma$ on $\R^N$, where $\gamma$ is the standard Gaussian density ($\gamma_k$ on $\R^k$), $0<a\le m\le b<\infty$, $\|\nabla m\|_\infty<\infty$, and $m$ together with its prefix averages
\[
 m_t(z_{1:t})=\int m(z_{1:t},z_{t+1:T})\gamma_{d(T-t)}(z_{t+1:T})\,dz_{t+1:T}
 \qquad (m_0=1,\ m_T=m)
\]
are smooth. Indeed, the conditionals $k_t(z_{1:t})=\gamma_d(z_t)m_t(z_{1:t})/m_{t-1}(z_{1:t-1})$ inherit \eqref{gpgeo:conditionalpoincare} from the Gaussian $L^q$-Poincar\'e inequality with $C_q=2(b/a)^{1/q}(\pi/2)(\E|Z|^q)^{1/q}$ ($Z\sim\NN(0,1)$), while the history and future scores
\[
 \nabla_i\log k_t=\nabla_i\log m_t-\nabla_i\log m_{t-1},
 \qquad
 \nabla_t\log(\widetilde\rho/\widetilde\rho_{1:t})=\nabla_t\log m-\nabla_t\log m_t
\]
are bounded in magnitude by $2\|\nabla m\|_\infty/a$.
\end{remark}

\subsection{Two-endpoint tangent expansions}

We now prove Proposition~\ref{prop:twopoint}, beginning with the general lower bound and then constructing matching bicausal upper bounds in the compact, log-concave reference, and scalar exponential-moment regimes.

\begin{lemma}[Local lower bound]
\label{gpgeo:lowerbound}
Let $P$ have density $\rho>0$, and let $h,k$ have zero integral with $h/\rho,k/\rho\in L^\infty$ and $N_{P,p}(k-h)<\infty$. For $P_\varepsilon^h=(\rho+\varepsilon h)\,dx$ and $P_\varepsilon^k=(\rho+\varepsilon k)\,dx$,
\[
 \liminf_{\varepsilon\to0}
 \frac{\AW_p(P_\varepsilon^h,P_\varepsilon^k)}{|\varepsilon|}
 \ge N_{P,p}(k-h).
\]
\end{lemma}
\begin{proof}
When the $\liminf$ is finite, take a subsequence achieving it, pick bicausal couplings $\pi_\varepsilon\in\Pi_{\rm bc}(P_\varepsilon^h,P_\varepsilon^k)$ of $(X,Y)$ whose costs exceed $\AW_p(P_\varepsilon^h,P_\varepsilon^k)^p$ by at most $o(|\varepsilon|^p)$, and set
\[
 V_t^\varepsilon(x)=\E_{\pi_\varepsilon}[(Y_t-X_t)/\varepsilon\mid X=x].
\]
Causality of $\pi_\varepsilon$ ensures that $V_t^\varepsilon(x)$ depends only on $x_{1:t}$, while $(\rho+\varepsilon h)/\rho=1+O(|\varepsilon|)$ and Jensen's inequality keep $V^\varepsilon$ bounded in $\mathcal A_{P,p}$. Passing to a subsequence, we obtain a weak limit $V^\varepsilon\rightharpoonup V$ in $\mathcal A_{P,p}$. For any $\phi\in C_c^\infty(\R^N)$, Taylor's theorem gives
\[
 |\phi(y)-\phi(x)-\nabla\phi(x)\cdot(y-x)|
 \le C_\phi\min\{|y-x|^2,|y-x|\}.
\]
Since $\int|y-x|^p\,d\pi_\varepsilon=O(|\varepsilon|^p)$, using $\min\{r^2,r\}\le r^p$ when $1<p<2$ and H\"older's inequality when $p\ge 2$ shows that
\[
 \frac1{|\varepsilon|}\int\min\{|y-x|^2,|y-x|\}\,d\pi_\varepsilon
 \le C\begin{cases}
 |\varepsilon|^{p-1},&1<p<2,\\
 |\varepsilon|,&p\ge2,
 \end{cases}
 \longrightarrow0.
\]
Dividing $\E_{\pi_\varepsilon}[\phi(Y)-\phi(X)]$ by $\varepsilon$ and letting $\varepsilon\to 0$ therefore gives
\[
 \int\phi(k-h)
 =\lim_{\varepsilon\to0}\int\sum_{t=1}^T\nabla_t\phi\cdot V_t^\varepsilon(\rho+\varepsilon h)
 =\int\sum_{t=1}^T\nabla_t\phi\cdot V_t\rho.
\]
Hence $-\sum_{t=1}^T\nabla_t\cdot(\rho V_t)=k-h$, and weak lower semicontinuity of $\|\cdot\|_{P,p}$ gives the lower bound.
\end{proof}

\paragraph{Compactly supported laws.}\label{sec:compact}

We first verify the conditional Poincar\'e, score, and comparison estimates for Gaussian mixtures with compactly supported mixing laws.

\begin{lemma}[Uniform mixture Poincar\'e inequalities]
\label{gpgeo:mixturepoincare}
If $\nu$ is a mixture of $\NN(a,\sigma^2 I_d)$ supported on location parameters $a$ of diameter at most $D$, then for $1\le q<\infty$ and $Z\sim\NN(0,1)$,
\begin{equation}
 \|\phi-\E_\nu\phi\|_{L^q(\nu)}
 \le \left(\frac\pi2(\E|Z|^q)^{1/q}\sigma
            +D e^{D^2/(8q\sigma^2)}\right)
                \|\nabla\phi\|_{L^q(\nu)}.
 \label{gpgeo:mixtureconstant}
\end{equation}
\end{lemma}
\begin{proof}
Let $B$ follow the mixing law independently of $Z_d\sim\NN(0,I_d)$, and put $m(a)=\E\phi(a+\sigma Z_d)$. Pisier's rotation interpolation on $[0,\pi/2]$ gives
\[
 \|\phi(a+\sigma Z_d)-m(a)\|_{L^q}
 \le \frac\pi2(\E|Z|^q)^{1/q}\sigma\|\nabla\phi(a+\sigma Z_d)\|_{L^q}.
\]
Let $B'$ be an independent copy of $B$. For any two locations $a,a'$ with $|a-a'|\le D$ and any $s\in[0,1]$, completing the square in the Gaussian density gives
\[
 \gamma_\sigma(x-[(1-s)a+sa'])
 \le e^{D^2/(8\sigma^2)}
       [(1-s)\gamma_\sigma(x-a)+s\gamma_\sigma(x-a')].
\]
Integrating $\nabla m$ along $(1-s)B+sB'$ therefore yields
\[
 \E|m(B)-m(B')|^q\le D^q e^{D^2/(8\sigma^2)}\int|\nabla\phi|^q\,d\nu,
\]
and the triangle inequality in $L^q(\nu)$ gives \eqref{gpgeo:mixtureconstant}.
\end{proof}

Whenever $\mu$ is supported in a compact set $K\subset\R^{dT}$ and $P=\mu*\gamma_\sigma$, Lemma~\ref{gpgeo:mixturepoincare} provides a conditional Poincar\'e constant $C_q$ in \eqref{gpgeo:conditionalpoincare} that is uniform in $x_{1:t-1}$. Differentiating the Gaussian kernel under the integral also yields the conditional score identities
\begin{align}
 \nabla_i\log k_t(x_{1:t})
 &=\sigma^{-2}\bigl(\E[A_i\mid X_{1:t}]-\E[A_i\mid X_{1:t-1}]\bigr),\quad i<t,
 \notag\\
 \nabla_t\log\frac{\rho(x)}{\rho_{1:t}(x_{1:t})}
 &=\sigma^{-2}\bigl(\E[A_t\mid X_{1:T}]-\E[A_t\mid X_{1:t}]\bigr),
 \label{gpgeo:futurescore}
\end{align}
both of which are bounded in $L^\infty$ by $\sigma^{-2}\operatorname{diam}(K)$.

\begin{proposition}[Comparison on a fixed compact latent class]
\label{gpgeo:compactcomparison}
Suppose $P=\nu*\gamma_\sigma$ and $Q=\widetilde\nu*\gamma_\sigma$ have densities $\rho$ and $\eta$, where $\nu$ and $\widetilde\nu$ are supported in the same compact set $K\subset\R^{dT}$. Then
\[
 \AW_p(P,Q)
 \le C(p,d,T,K,\sigma)\left\|\frac{\eta-\rho}{\rho}\right\|_{L^p(P)}.
\]
\end{proposition}
\begin{proof}
Set $u=\eta/\rho$ and $u_t=\E_P[u\mid X_{1:t}]$ as in Lemma~\ref{gpgeo:onesided}, and split the history space at date $t-1$ into $\{u_{t-1}\ge 1/2\}$ and $\{u_{t-1}<1/2\}$. On $\{u_{t-1}\ge 1/2\}$, Lemma~\ref{gpgeo:ordinarycomparison} bounds the conditional $W_p^p$ integral by $(pC_q)^p 2^{p-1}\E_P|u_t-u_{t-1}|^p\le C\|u-1\|_{L^p(P)}^p$. On $\{u_{t-1}<1/2\}$, both $k_t^P(z,\cdot)$ and $k_t^Q(z,\cdot)$ are mixtures of $\NN(a,\sigma^2 I_d)$ with $a$ in the $t$-th projection $K_t$, so $W_p(k_t^P(z,\cdot),k_t^Q(z,\cdot))\le D_t:=\operatorname{diam}(K_t)$. Since $u_{t-1}<1/2$ implies $|u_{t-1}-1|>1/2$, we obtain
\begin{align*}
 D_t^p\E_P[u_{t-1}\mathbf1_{\{u_{t-1}<1/2\}}]
 &\le D_t^p P_{1:t-1}(u_{t-1}<1/2)\\
 &\le 2^pD_t^p\E_P|u_{t-1}-1|^p
 \le 2^pD_t^p\|u-1\|_{L^p(P)}^p.
\end{align*}
Inserting both bounds into the coupling recursion of Lemma~\ref{gpgeo:onesided}, together with \eqref{gpgeo:historylipschitz} and \eqref{gpgeo:futurescore}, gives the stated inequality.
\end{proof}

\begin{proposition}[Two-endpoint tangent expansion under compact support]
\label{gpgeo:compactexpansion}
Under Assumption~\ref{ass:compact}, for all $h,k\in\DD$ as in \eqref{eq:tilts},
\[
 \lim_{\varepsilon\to0}
 \frac{\AW_p(P_\varepsilon^h,P_\varepsilon^k)}{|\varepsilon|}
       =N_{P,p}(k-h).
\]
\end{proposition}
\begin{proof}
Let $v\in\mathcal A_{P,p}$ be the minimal-norm solution to $D_\rho v=(k-h)/\rho$ from Lemma~\ref{gpgeo:rightinverse}, and let $w=(w_1,\dots,w_T)$ with $w_t\in C_c^\infty(\R^{dt};\R^d)$ be a smooth approximant from Lemma~\ref{gpgeo:core}. The ODE flow $\partial_\varepsilon F_\varepsilon(x)=w(F_\varepsilon(x))$ with $F_0(x)=x$ is a global block-triangular diffeomorphism with block-triangular inverse $F_{-\varepsilon}$, so $(\operatorname{Id},F_\varepsilon)_\# P_\varepsilon^h$ is a bicausal coupling between $P_\varepsilon^h$ and $Q_\varepsilon^w:=(F_\varepsilon)_\# P_\varepsilon^h$. Since $(F_\varepsilon-\operatorname{Id})/\varepsilon\to w$ uniformly, dominated convergence gives
\[
 \frac1{|\varepsilon|^p}
 \int\sum_{t=1}^T|F_{\varepsilon,t}(x)-x_t|^p(\rho+\varepsilon h)(x)\,dx
 \longrightarrow\|w\|_{P,p}^p.
\]
Put $a=h/\rho$ and $b_w=D_\rho w$. Because $\nabla_t\log\rho=\nabla_t\log\rho_{1:t}+\nabla_t\log(\rho/\rho_{1:t})$ is bounded on the compact support of $w_t$ by \eqref{gpgeo:futurescore}, $b_w$ is bounded and continuous on $\R^N$, and the change-of-variables formula along $F_\varepsilon$ yields
\[
 \frac{dQ_\varepsilon^w}{dP}(x)
 =[1+\varepsilon a(F_{-\varepsilon}(x))]
     \exp\!\left(\int_0^\varepsilon b_w(F_{-s}(x))\,ds\right).
\]
Thus $dQ_\varepsilon^w/dP=1+O(|\varepsilon|)$ uniformly on $\R^N$ and
\[
 \frac{dQ_\varepsilon^w/dP-1}{\varepsilon}
          \longrightarrow a+b_w\quad\text{in }L^p(P).
\]
For small $|\varepsilon|$, we have $dQ_\varepsilon^w/dP_\varepsilon^k\ge 1/2$. Applying Lemma~\ref{gpgeo:onesided} with $P_\varepsilon^k$ as the reference law and passing to the limit as $\varepsilon\to 0$ gives
\[
 \limsup_{\varepsilon\to0}
 \frac{\AW_p(Q_\varepsilon^w,P_\varepsilon^k)}{|\varepsilon|}
 \le C\left\|D_\rho w-\frac{k-h}{\rho}\right\|_{L^p(P)},
\]
with $C$ independent of $w$. Combining this with the transport cost of $F_\varepsilon$, letting $w\to v$ in the graph norm of Lemma~\ref{gpgeo:core}, and matching with Lemma~\ref{gpgeo:lowerbound} proves the expansion.
\end{proof}

\begin{corollary}[A common Hilbert control on compact latent classes]
\label{gpgeo:hilbertcomparison}
Fix an integer $s>N/2$, and choose $0<a<(2\sigma^2)^{-1}<\alpha$ with $pa>(p-1)\alpha$. On the Hilbert space $\HH_{a,s}=\{h:e^{a|\cdot|^2}h\in H^s(\R^N)\}$ with norm $\|h\|_{\HH_{a,s}}=\|e^{a|\cdot|^2}h\|_{H^s(\R^N)}$ and $\HH_{a,s,0}=\{h\in\HH_{a,s}:\int h=0\}$, all $P=\nu*\gamma_\sigma$ and $Q=\widetilde\nu*\gamma_\sigma$ with $\supp(\nu),\supp(\widetilde\nu)\subset K$ compact satisfy
\[
 \AW_p(P,Q)\le C\|\rho-\eta\|_{\HH_{a,s}},
 \qquad
 N_{P,p}(h)\le C\|h\|_{\HH_{a,s}}\quad\Bigl(\int h=0\Bigr).
\]
\end{corollary}
\begin{proof}
Since $\nu$ is supported in the compact set $K$, we have $\rho(x)\ge c_\alpha e^{-\alpha|x|^2}$, while Sobolev embedding ($s>N/2$) gives $|h(x)|\le C_s\|h\|_{\HH_{a,s}}e^{-a|x|^2}$. Hence
\[
 \int |h|^p\rho^{1-p}
 \le C\|h\|_{\HH_{a,s}}^p
 \int e^{-[pa-(p-1)\alpha]|x|^2}\,dx
 \le C'\|h\|_{\HH_{a,s}}^p,
\]
and the two bounds follow from Proposition~\ref{gpgeo:compactcomparison} and Lemma~\ref{gpgeo:rightinverse}.
\end{proof}

\paragraph{Strongly log-concave reference laws.}\label{sec:reference}

When $\mu=f\mu_0$ satisfies Assumption~\ref{ass:structured}, we construct a global bi-Lipschitz triangular diffeomorphism $F$ pushing the standard Gaussian $\gamma$ on $\R^N$ onto $P_0=\mu_0*\gamma_\sigma$, pull the adapted continuity equation back to Gaussian chart coordinates, and conjugate the ODE recovery flow through $F$.

Write $\rho_0=dP_0/dx$, enumerate the $N$ scalar coordinates in lexicographic order $(t,i)$, and let $\bar\rho_{0,j}(x_{1:j})$ be the $j$-th scalar prefix density of $P_0$ ($\bar\rho_{0,0}=1$). Under $A\sim\mu_0$ and $X=A+\sigma Z$, the conditional law of $A$ given $X_{1:j}=x\in\R^j$ has potential
\[
 W_{j,x}(a)=V_0(a)+\frac{|x-a_{1:j}|^2}{2\sigma^2}+\text{const},
 \qquad
 \nabla^2W_{j,x}\succeq K_j:=\kappa I_N+\sigma^{-2}\Pi_j,
\]
where $\Pi_j$ projects onto the first $j$ scalar coordinates. By the Brascamp--Lieb variance inequality \cite[Eq.~(1.3)]{CCL}, $\Cov(A\mid X_{1:j}=x)\preceq K_j^{-1}$, so $\Cov(A_{1:j}\mid X_{1:j}=x)\preceq(\kappa+\sigma^{-2})^{-1}I_j$. Differentiating $\bar\rho_{0,j}(x)=\E_{\mu_0}[\gamma_\sigma(x-A_{1:j})]$ twice gives $\nabla^2(-\log\bar\rho_{0,j})(x)=\sigma^{-2}I_j-\sigma^{-4}\Cov(A_{1:j}\mid X_{1:j}=x)$, and therefore
\begin{equation}
 mI_j\preceq\nabla^2(-\log\bar \rho_{0,j})\preceq MI_j,
 \qquad
 m=\frac{\kappa}{1+\kappa\sigma^2},\quad M=\sigma^{-2}.
 \label{eq:Hessianbounds}
\end{equation}

For $1\le j\le N$ and a scalar history $a\in\R^{j-1}$, the conditional density $k_a(z)=\bar\rho_{0,j}(a,z)/\bar\rho_{0,j-1}(a)$ and potential $U_a=-\log k_a$ satisfy
\[
 m\le U_a''\le M,
 \qquad |U_a'(z)-U_b'(z)|\le M|a-b|.
\]
With $r=M|a-b|/m$, the inequality $U_b'(z)\ge U_a'(z-r)$ implies that $k_b(z)/k_a(z-r)$ is nonincreasing in $z$, so the conditional quantile $q_a(u)=\inf\{z:\int_{-\infty}^z k_a(y)\,dy\ge u\}$ obeys
\[
 |q_a(u)-q_b(u)|\le B|a-b|,
 \qquad B=M/m,\quad 0<u<1.
\]
Moreover, Caffarelli's contraction theorem \cite{Caf,FGP} applied to $T_a(z)=q_a(\Phi(z))$ and $T_a^{-1}$ gives $M^{-1/2}\le\partial_zT_a(z)\le m^{-1/2}$. Defining $F_1(z_1)=T_\varnothing(z_1)$ and $F_j(z_{1:j})=T_{F_{1:j-1}(z)}(z_j)$ produces a smooth triangular diffeomorphism $F:\R^N\to\R^N$ with $F_\#\gamma=P_0$ and inverse $(F^{-1})_j(x)=T_{x_{1:j-1}}^{-1}(x_j)$. Writing $\ell=m^{-1/2}$, the quantile bounds yield
\[
 |F_j(z)-F_j(z')|
 \le B|F_{1:j-1}(z)-F_{1:j-1}(z')|+\ell|z_j-z_j'|,
\]
whence $\operatorname{Lip}(F)^2\le 2m^{-1}(1+2B^2)^{N-1}$, $\operatorname{Lip}(F^{-1})^2\le M\{N+(N-1)B^2\}$, and consequently $\sup_z\|DF(z)\|<\infty$ and $\sup_x\|D(F^{-1})(x)\|<\infty$.

For $\mu=f\mu_0$, we have $\rho(x)=\rho_0(x)r(x)$ with $r(x)=\E_{\mu_0}[f(A)\mid A+\sigma Z=x]\in[c_f,C_f]$. For any $g\in L^\infty(\mu_0)$ and unit vector $e\in\R^N$, the variance bound above gives
\[
 \left|\partial_e\E_{\mu_0}[g(A)\mid X=x]\right|
 = \sigma^{-2}\big|\Cov_{\mu_0}(g(A),e\cdot A\mid X=x)\big|
 \le\frac{\|g\|_\infty}{\sigma^2\sqrt{\kappa+\sigma^{-2}}}.
\]
In the chart $F$, the law $P$ and its tilts $P_\varepsilon^h$ ($h=(g\mu)*\gamma_\sigma\in\DD$) take the form
\[
 P=F_\#\widetilde P,\qquad \widetilde\rho(z)=\gamma(z)\widetilde r(z),\qquad \widetilde r=r\circ F,
\]
and $P_\varepsilon^h=F_\#(\widetilde\rho_\varepsilon^h\,dz)$ with $\widetilde\rho_\varepsilon^h=\gamma(\widetilde r+\varepsilon\widetilde r_h)$ and $\widetilde r_h(z)=\E_{\mu_0}[f(A)g(A)\mid X=F(z)]$. Both $\widetilde r$ and $\widetilde r_h$ are bounded with bounded gradients, and their block-prefix averages are smooth as well—for $\widetilde r$, the $t$-th block-prefix average equals $(\rho_{1:t}/\rho_{0,1:t})\circ F_{1:t}$, with the same identity for the tilted densities. Remark~\ref{gpgeo:gaussianmultiplier} therefore applies uniformly to $\widetilde\rho_\varepsilon^h$ for small $|\varepsilon|$.

Given a block-adapted field $v\in\mathcal A_{P,p}$, define the pulled-back velocity $u(z)=DF(z)^{-1}v(F(z))$ and density direction $\widetilde h(z)=h(F(z))|\det DF(z)|$. Since $DF$ and $DF^{-1}$ are lower block triangular, $u$ is block adapted and
\[
 h=-\diver(\rho v)\quad\Longleftrightarrow\quad \widetilde h=-\diver(\widetilde\rho u),
\]
while the physical kinetic energy becomes
\begin{equation}
 \|v\|_{P,p}^p
 =\int\sum_{t=1}^T|(DF(z)u(z))_t|^p\,\widetilde\rho(z)\,dz.
 \label{eq:variablemetric}
\end{equation}
By the uniform bounds on $\|DF\|$ and $\|DF^{-1}\|$, the $p$-th root of \eqref{eq:variablemetric} defines an equivalent strictly convex norm on $\mathcal A_{\widetilde P,p}$. Applying Remark~\ref{gpgeo:gaussianmultiplier} and Lemmas~\ref{gpgeo:rightinverse}--\ref{gpgeo:core} to $\widetilde\rho=\gamma\widetilde r$, together with $\int|\widetilde h/\widetilde\rho|^p\widetilde\rho=\int|h|^p\rho^{1-p}$, yields
\begin{equation}
 N_{P,p}(h)\le C_{p,P}\|h/\rho\|_{L^p(P)},
 \label{eq:referenceLpbound}
\end{equation}
along with a unique minimizer $u$ of \eqref{eq:variablemetric} and smooth block-prefix-compact approximants $u^m\to u$ in $\mathcal A_{\widetilde P,p}$ satisfying $D_{\widetilde\rho}u^m\to\widetilde h/\widetilde\rho$ in $L^p(\widetilde P)$.

Now let $h,j\in\DD$, let $u$ be the pulled-back minimal velocity for $j-h$, let $\Phi_\varepsilon^m$ be the ODE flow of $u^m$, and put $\Psi_\varepsilon^m=F\circ\Phi_\varepsilon^m\circ F^{-1}$. Because both $F$ and $\Phi_\varepsilon^m$ are block triangular with block-triangular inverses, $(\operatorname{Id},\Psi_\varepsilon^m)_\# P_\varepsilon^h$ is bicausal, and dominated convergence gives
\[
 |\varepsilon|^{-p}\int\sum_{t=1}^T|\Psi_{\varepsilon,t}^m(x)-x_t|^p\,dP_\varepsilon^h(x)
 \longrightarrow
 \int\sum_{t=1}^T|(DF(z)u^m(z))_t|^p\,d\widetilde P(z).
\]
In chart coordinates, the transported density $\widetilde q_\varepsilon^m=d((\Phi_\varepsilon^m)_\#\widetilde P_\varepsilon^h)/dz$ satisfies
\[
 \frac{\widetilde q_\varepsilon^m-\widetilde\rho}{\varepsilon}
 \longrightarrow \widetilde h+\widetilde\rho D_{\widetilde\rho}u^m
 \quad\text{in }L^p(\widetilde\rho^{1-p}\,dz).
\]
Applying Lemma~\ref{gpgeo:onesided} with reference density $\widetilde\rho+\varepsilon\widetilde j$ and pushing the resulting bicausal coupling forward by $F$ (which multiplies $\AW_p$ by at most $L_p=T^{|1/p-1/2|}\operatorname{Lip}(F)$) yields
\[
 \limsup_{\varepsilon\to 0}
 \frac{\AW_p((\Psi_\varepsilon^m)_\# P_\varepsilon^h,P_\varepsilon^j)}{|\varepsilon|}
 \le C L_p\left\|D_{\widetilde\rho}u^m-
 \frac{\widetilde j-\widetilde h}{\widetilde\rho}\right\|_{L^p(\widetilde P)}.
\]
Letting $\varepsilon\to 0$ first and then $m\to\infty$, and combining with Lemma~\ref{gpgeo:lowerbound}, establishes \eqref{eq:twopoint} under Assumption~\ref{ass:structured}.

\paragraph{Scalar paths with exponential moments.}\label{gpstat:sec-scalar-low}

Suppose $d=1$ and the assumptions of Theorem~\ref{thm:scalar} hold, and set $r=\max(p,2)$. By the weighted energy estimates in Section~\ref{sec:optimizedtails} (\eqref{eq:energybound} when $1<p\le 2$, $r=2$) and Section~\ref{gpstat:sec-high-p} (\eqref{gpstat:eq-energy} when $p>2$, $r=p$), the recursive scalar velocity $v_h$ in \eqref{eq:scalarrecursion} satisfies $\|v_h\|_{P,p}\le T^{(1/p-1/2)_+}\|v_h\|_{P,r}<\infty$. In one dimension, $v_h$ is the unique field in $\mathcal A_{P,p}$ solving $-\sum_{t=1}^T\partial_t(\rho v_t)=h$: all fluxes $\rho v_t$ lie in $L^1(\R^T)$, so integrating out future coordinates shows by induction on $1\le j\le T$ that the difference of any two solutions has a flux independent of $x_j$, which must vanish by $L^1(\R)$ integrability. Hence $Q_{P,p}(h)=\sum_{t=1}^T\int|v_{h,t}|^p\rho_{1:t}$.

For a positive smooth density $f$ on $\R^T$, let $R^f:(0,1)^T\to\R^T$ denote its inverse Rosenblatt map, defined by $R_1^f(u)=F_{f,1}^{-1}(u_1)$ and $R_t^f(u)=F_{f,t\mid R_{1:t-1}^f(u)}^{-1}(u_t)$. Given $h=h_g\in\DD$, set $\rho_\varepsilon=\rho+\varepsilon h$ and $R_\varepsilon=R^{\rho_\varepsilon}$, and choose $\varepsilon_0>0$ so that $\varepsilon_0\|g\|_\infty\le 1/2$ and all $\rho_\varepsilon$ ($|\varepsilon|\le\varepsilon_0$) lie in a single moment class $\CC_{M_0}$ (defined in Section~\ref{sec:optimizedtails}). Differentiating the conditional quantile equations yields
\[
 w_{\varepsilon,t}(u):=\partial_\varepsilon R_{\varepsilon,t}(u)
 =v_{\rho_\varepsilon,h,t}(R_\varepsilon(u)_{1:t})
 \longrightarrow v_{\rho,h,t}(R_0(u)_{1:t})
\]
pointwise on $(0,1)^T$ as $\varepsilon\to 0$.

By Gaussian convolution, positivity of the prefix densities, and the recursion \eqref{eq:scalarrecursion}, each velocity $v_{\rho_\varepsilon,h,t}(x_{1:t})$ is jointly continuous in $(\varepsilon,x_{1:t})$, so
\[
 M_R:=\max_{1\le t\le T}
 \sup_{\substack{|\varepsilon|\le\varepsilon_0\\|x_{1:t}|\le R}}
 |v_{\rho_\varepsilon,h,t}(x_{1:t})|<\infty
 \qquad\text{for every }R<\infty.
\]
For $M>M_R$, changing variables via the prefix Rosenblatt map and applying \eqref{eq:energybound} (when $r=2$) or \eqref{gpstat:eq-energy} (when $r=p>2$) gives
\begin{equation}
\begin{split}
 &\int_{(0,1)^T}|w_{\varepsilon,t}(u)|^r
   \mathbf1_{\{|w_{\varepsilon,t}(u)|>M\}}\,du\\
 &\quad\le\int_{|x_{1:t}|>R}
   |v_{\rho_\varepsilon,h,t}(x_{1:t})|^r
   (\rho_\varepsilon)_{1:t}(x_{1:t})\,dx_{1:t}
 \le C e^{-\eta R^2}\|h\|_\HH^r,
\end{split}
\label{eq:scalar-uniform-integrability}
\end{equation}
with constants uniform in $|\varepsilon|\le\varepsilon_0$. Choosing $R$ large and then $M>M_R$ shows that $\{|w_{\varepsilon,t}|^r\}_{|\varepsilon|\le\varepsilon_0}$ is uniformly integrable, so Vitali's theorem yields $w_{\varepsilon,t}\to v_{\rho,h,t}\circ R_0$ in $L^r((0,1)^T)$. Integrating in $\varepsilon$ gives
\[
 \frac{R_\varepsilon-R_0}{\varepsilon}
 =\int_0^1\left.\partial_\delta R_\delta\right|_{\delta=s\varepsilon}\,ds
 \longrightarrow v_{\rho,h}\circ R_0
 \quad\text{in }L^r((0,1)^T;\ell^r),
\]
and hence also in $L^p((0,1)^T;\ell^p)$ since $p\le r$. Applying this expansion to $h$ and $k$ and using the bicausal coupling $(R_\varepsilon^h(U),R_\varepsilon^k(U))$ with $U\sim\mathrm{Unif}((0,1)^T)$ yields the upper bound in \eqref{eq:twopoint}; Lemma~\ref{gpgeo:lowerbound} supplies the matching lower bound.

\subsection{Statistical limits and bootstrap}\label{sec:densityspaces}\label{gpstat:sec-statistics}\label{sec:statistics}

We specify the auxiliary Hilbert space $\HH$ in each of the four regimes so that $\E_\mu\|e_A\|_\HH^2<\infty$, $N_{P,p}$ is continuous on $\HH_0$, and $\AW_p$ is controlled by $\|\cdot\|_\HH$:
\begin{itemize}
\item \emph{Compact support (Assumption~\ref{ass:compact}, $1<p<\infty$).} Fix an integer $s>N/2$, pick $0<a<(2\sigma^2)^{-1}<\alpha$ with $pa>(p-1)\alpha$, and set $\HH=\HH_{a,s}$ from Corollary~\ref{gpgeo:hilbertcomparison}. Then $\sup_{y\in[-R,R]^N}\|e_y\|_{\HH_{a,s}}<\infty$, while Corollary~\ref{gpgeo:hilbertcomparison} yields $\AW_p(\rho\,dx,\eta\,dx)\le C\|\eta-\rho\|_\HH$ and $N_{P,p}(h)\le C\|h\|_\HH$ uniformly over mixing laws supported in $[-R,R]^N$.

\item \emph{Scalar exponential moment (Assumption~\ref{ass:scalar}, $1<p\le 2$, $\beta\sigma^2>\theta_T$).} Take $\HH=\HH_\omega=L^2(\R^T,e^{\omega|x|^2}\,dx)$ from Section~\ref{sec:optimizedtails} with
\[
 \frac1{2\sigma^2}+\frac{T-1}{\sigma^2(1+2\beta\sigma^2)}
 <\omega<\frac{\beta}{1+\beta\sigma^2}.
\]
Direct integration gives
\[
 \|e_y\|_{\HH_\omega}^2
 =[4\pi\sigma^2(1-\omega\sigma^2)]^{-T/2}
 \exp\!\left(\frac{\omega|y|^2}{1-\omega\sigma^2}\right),
\]
where $\omega/(1-\omega\sigma^2)<\beta$, so $\E_\mu\|e_A\|_{\HH_\omega}^2<\infty$. On any moment class $\CC_{M_0}$, Proposition~\ref{prop:metric} together with H\"older's inequality
\begin{equation}
 \AW_p\le T^{1/p-1/2}\AW_2,
 \qquad N_{P,p}\le T^{1/p-1/2}N_{P,2}
 \label{gpstat:eq-low-p-domination}
\end{equation}
provides the metric comparison and continuity of $N_{P,p}$.

\item \emph{Log-concave reference (Assumption~\ref{ass:structured}, $1<p\le 2$).} Take $\HH=\HH_P=L^2(\R^N,\rho^{-1}\,dx)$ and $\HH_{P,0}=\{h\in\HH_P:\int h=0\}$. Lemma~\ref{lem:intrinsicmoment} gives $\E_\mu\|e_A\|_{\HH_P}^2<\infty$ at every $\sigma>0$, while \eqref{eq:regularnetcomparison}, \eqref{eq:referenceLpbound}, and \eqref{gpstat:eq-low-p-domination} give the metric comparison and continuity of $N_{P,p}$.

\item \emph{Unbounded support with $p>2$ (Assumption~\ref{ass:scalar} or~\ref{ass:structured}, $\beta\sigma^2>\Theta_{p,N}$).} Take $\HH=\HH_{a,s}$ with $(\alpha,b,a)$ satisfying \eqref{gpstat:eq-parameters}; the feature moment and metric comparisons are given by \eqref{gpstat:eq-high-second-moment} and Proposition~\ref{gpstat:prop-comparison}.
\end{itemize}

In all four regimes, the centered kernel features $\xi_i=e_{A^{(i)}}-\rho$ are i.i.d.\ square-integrable random elements of $\HH_0$, so the Hilbert-space central limit theorem \cite{LT} yields
\begin{equation}
 G_n:=\sqrt n(\rho_n-\rho)
 =\frac1{\sqrt n}\sum_{i=1}^n\xi_i
 \Longrightarrow G\quad\text{in }\HH,
 \label{gpstat:eq-hilbert-clt}
\end{equation}
where $G$ is centered Gaussian with covariance $C_Gu=\E[\langle\xi_1,u\rangle_\HH\xi_1]$.

\paragraph{Uniform local expansion and weak limit.}

Define $\mathsf T:L^2_0(\mu)\to\HH$ on $L^2_0(\mu)=\{g\in L^2(\mu):\int g\,d\mu=0\}$ by $\mathsf Tg=\E[g(A)(e_A-\rho)]$. Because $C_G=\mathsf T\mathsf T^*$ and bounded zero-mean functions are dense in $L^2_0(\mu)$, the topological support of $G$ in $\HH$ is
\[
 \mathcal S=\overline{\operatorname{Ran}\mathsf T}^{\,\HH}
           =\overline{\DD}^{\,\HH},
\]
and $\xi_i\in\mathcal S$ almost surely. In the exponential-moment regimes, where the metric comparisons apply on a bounded moment class $\CC_{M_0}$, put $M_i=e^{\beta|A^{(i)}|^2}$, $m_\beta=\E M_1<M_0$, and $M_i^K=\min(M_i,K)$. Then $n^{-1}\sum_{i=1}^n M_i\to m_\beta$ in probability and
\[
 \Var^*\left(\frac1n\sum_{j=1}^nM_{I_j}^K\right)
 \le\frac{K^2}{n},\qquad
 \E\E^*\left[\frac1n\sum_{j=1}^n(M_{I_j}-M_{I_j}^K)\right]
 =\E(M_1-M_1^K),
\]
so sending $n\to\infty$ and then $K\to\infty$ gives $\Pp(P_n\in\CC_{M_0},\,P_n^*\in\CC_{M_0})\to 1$; moreover, $m_\beta+\varepsilon\E_\mu[g(A)e^{\beta|A|^2}]<M_0$ for any $g\in L^\infty(\mu)$ once $|\varepsilon|$ is small.

For directions $a,b\in\HH_0$ such that $\rho+\varepsilon a$ and $\rho+\varepsilon b$ are feasible probability densities (in the compact class, $\CC_{M_0}$, or $\HH_P$), set $F_\varepsilon(a,b)=\varepsilon^{-1}\AW_p((\rho+\varepsilon a)\,dx,(\rho+\varepsilon b)\,dx)$. We claim that for every compact $K\subset\mathcal S^2$,
\begin{equation}
 \sup_{\substack{(a,b)\in K\\(a,b)\text{ feasible at }\varepsilon}}
 |F_\varepsilon(a,b)-N_{P,p}(b-a)|\longrightarrow0\qquad (\varepsilon\downarrow 0).
 \label{gpstat:eq-compact-uniform}
\end{equation}
Indeed, fix $\delta>0$ and, by density of $\DD$ in $\mathcal S$, pick a finite $\delta$-net $\{(a_\ell,b_\ell)\}_{\ell=1}^M\subset\DD^2$ for $K$. For $\varepsilon>0$ small enough, all $\rho+\varepsilon a_\ell$ and $\rho+\varepsilon b_\ell$ are feasible, and the metric comparisons established above give
\[
 |F_\varepsilon(a,b)-F_\varepsilon(a_\ell,b_\ell)|
 \le L_{\mathrm{cmp}}(\|a-a_\ell\|_{\HH}+\|b-b_\ell\|_{\HH}),
\]
with $L_{\mathrm{cmp}}$ independent of the net and of $\varepsilon$: in the compact and exponential-moment regimes all endpoints belong to one fixed compact class or to $\CC_{M_0}$, while in the reference regime ($1<p\le 2$) choosing $|\varepsilon|\max_g\|g\|_\infty\le 1/2$ over the finitely many tilts $h_g$ in the net keeps the tilted multipliers in $[c_f/2,3C_f/2]$, so \eqref{eq:regularnetcomparison} supplies a uniform constant. Writing $L_{\mathrm{tan}}$ for the Lipschitz constant of $N_{P,p}$ on $\HH_0$, the triangle inequality gives
\begin{align*}
 &\sup_{\substack{(a,b)\in K\\\text{feasible at }\varepsilon}}
 \big|F_\varepsilon(a,b)-N_{P,p}(b-a)\big|\\
 &\quad\le \max_{1\le\ell\le M}
 \big|F_\varepsilon(a_\ell,b_\ell)-N_{P,p}(b_\ell-a_\ell)\big|
 +(L_{\mathrm{cmp}}+L_{\mathrm{tan}})\delta.
\end{align*}
For a fixed net, the maximum vanishes as $\varepsilon\downarrow 0$ by Proposition~\ref{prop:twopoint}; letting $\delta\downarrow 0$ proves \eqref{gpstat:eq-compact-uniform}. Since $\{G_n\}$ is tight in $\mathcal S$ by \eqref{gpstat:eq-hilbert-clt} and $P_n$ is feasible with probability tending to one, applying \eqref{gpstat:eq-compact-uniform} to $(a,b)=(0,G_n)$ with $\varepsilon=n^{-1/2}$ gives
\[
 \sqrt n\,\AW_p(P_n,P)-N_{P,p}(G_n)\longrightarrow0\quad\text{in probability},
\]
and continuous mapping yields \eqref{eq:mainlimit}.

\paragraph{Bootstrap, continuity, and moment convergence.}

Set $G_n^*=\sqrt n(\rho_n^*-\rho_n)$. We first check that in $\HH$,
\begin{equation}
 d_{\mathrm{BL},\HH}(\Law^*(G_n^*),\Law(G))\longrightarrow0
 \quad\text{in probability}.
 \label{gpstat:eq-hilbert-bootstrap}
\end{equation}
Let $\Pi_m$ be finite-rank orthogonal projections on $\HH$ increasing strongly to the identity. At each fixed $m$, $\E\|\xi_1\|_\HH^2<\infty$ yields the Lindeberg condition
\[
 \E\left[\frac1n\sum_{i=1}^n\|\xi_i\|_{\HH}^2
       \mathbf 1_{\{\|\xi_i\|_{\HH}>\delta\sqrt n\}}\right]
 \longrightarrow0,
\]
so the finite-dimensional bootstrap CLT \cite[Chapter~23]{vdV} applies to $\Pi_m G_n^*$. Meanwhile, the projection tail satisfies $\E^*\|(I-\Pi_m)G_n^*\|_\HH^2\le n^{-1}\sum_{i=1}^n\|(I-\Pi_m)\xi_i\|_\HH^2\to\E\|(I-\Pi_m)\xi_1\|_\HH^2\to 0$ as $n\to\infty$ and then $m\to\infty$, giving \eqref{gpstat:eq-hilbert-bootstrap}.

Because $(G_n,G_n+G_n^*)$ is tight in $\mathcal S^2$ and $(P_n,P_n^*)$ is feasible with probability tending to one, \eqref{gpstat:eq-compact-uniform} implies that
\[
 R_n:=\sqrt n\,\AW_p(P_n^*,P_n)-N_{P,p}(G_n^*)\longrightarrow0
\]
in joint probability. Since $\E[\Pp^*(|R_n|>\delta)]=\Pp(|R_n|>\delta)\to 0$, we have $\Pp^*(|R_n|>\delta)\to 0$ in probability for every $\delta>0$. For every $f:\R\to\R$ with $\|f\|_\infty\le 1$ and $\operatorname{Lip}(f)\le 1$,
\[
 \big|\E^*f(\sqrt n\,\AW_p(P_n^*,P_n))-\E^*f(N_{P,p}(G_n^*))\big|
 \le\delta+2\Pp^*(|R_n|>\delta),
\]
and combining this bound with \eqref{gpstat:eq-hilbert-bootstrap} and the continuity of $N_{P,p}$ proves \eqref{eq:bootstrap}.

To upgrade \eqref{eq:bootstrap} to Kolmogorov distance, we verify that $N_{P,p}(G)$ is atomless whenever $\mu$ is not a point mass. For any distinct $h,k$ in the domain of $Q_{P,p}=N_{P,p}^p$, their minimal adapted velocities $v_h$ and $v_k$ differ because $D_\rho v_h\neq D_\rho v_k$. Since $tv_h+(1-t)v_k$ is feasible for $th+(1-t)k$, strict convexity of $\|\cdot\|_{P,p}^p$ gives
\[
 Q_{P,p}(th+(1-t)k)
 \le \|tv_h+(1-t)v_k\|_{P,p}^p
 <tQ_{P,p}(h)+(1-t)Q_{P,p}(k),\qquad 0<t<1.
\]
When $\mu$ is not a point mass, $C_G\neq 0$. Choosing a unit eigenvector $e\in\HH_0$ of $C_G$ with eigenvalue $\lambda>0$ and setting $h=\sqrt\lambda e$, we decompose $G=Zh+G^\perp$ with $Z\sim\NN(0,1)$ independent of $G^\perp$. For each fixed $G^\perp$, the function $z\mapsto Q_{P,p}(G^\perp+zh)$ is strictly convex on $\R$, so every level set has at most two points and hence zero Gaussian measure. By Fubini's theorem, $N_{P,p}(G)$ has a continuous distribution function, and P\'olya's theorem upgrades \eqref{eq:bootstrap} to Kolmogorov distance.

It remains to verify the moment convergence \eqref{eq:compactmoments}. Under Assumption~\ref{ass:compact}, Corollary~\ref{gpgeo:hilbertcomparison} yields $\sqrt n\,\AW_p(P_n,P)\le C\|G_n\|_\HH$ samplewise. Since the features $\xi_i$ are bounded in $\HH$, the Rademacher moment inequality \cite{LT} gives $\sup_n\E\|G_n\|_\HH^r<\infty$ for all $r>\max(p,2)$, ensuring uniform integrability and \eqref{eq:compactmoments}. Under Assumption~\ref{ass:structured} with $1<p\le 2$, \eqref{eq:anchoredchisquare} and \eqref{gpstat:eq-low-p-domination} give the samplewise bound
\[
 \sqrt n\,\AW_p(P_n,P)
 \le T^{1/p-1/2}\sqrt{C_{\mathrm A}}\,\|G_n\|_{\HH_P}.
\]
Put $X_n=\|G_n\|_{\HH_P}^2$ and $X=\|G\|_{\HH_P}^2$. Then $X_n\Longrightarrow X$ and $\E X_n=\E\|\xi_1\|_{\HH_P}^2=\E X<\infty$ for every $n$, whence $\E(X_n-K)_+\to\E(X-K)_+$ for every $K>0$. Using $X_n\mathbf1_{\{X_n>2K\}}\le 2(X_n-K)_+$ establishes uniform integrability of $\{X_n\}$ and proves \eqref{eq:compactmoments} for $1<p\le 2$.

\paragraph{Proof of the Gaussian obstruction.}

\begin{proof}[Proof of Proposition~\ref{prop:failure}]
Let $A\sim\NN(0,s^2)$, $\tau^2=s^2+\sigma^2$, and $P=\NN(0,\tau^2)$, with density $\rho=\varphi_\tau$ and CDF $F(x)=\Phi(x/\tau)$, and write $F_n(x)=\int_{-\infty}^x\rho_n(y)\,dy$. Because $\E\int e_A^2/\rho=1+s^2/\sigma^2<\infty$, the Hilbert-space CLT in $\HH_{P,0}$ implies $\sqrt n(F_n-F)\Longrightarrow B_\sigma$ in $L^\infty(\R)$, where $B_\sigma(x)=\int_{-\infty}^x G(y)\,dy$ has covariance $\Cov(B_\sigma(x),B_\sigma(y))=\Cov_\mu(\Phi_\sigma(x-A),\Phi_\sigma(y-A))$. Writing $Y_i=A+\sigma Z_i$ ($Z_1,Z_2\sim\NN(0,1)$ i.i.d.\ and independent of $A$), we have $V(x):=\Var(B_\sigma(x))=\Pp(Y_1>x,Y_2>x)-\Pp(Y_1>x)^2$. Integrating the bivariate Gaussian density of $(Y_1,Y_2)$ over $[x,x+x^{-1}]^2$ for large $x>0$ gives $\Pp(Y_1>x,Y_2>x)\ge cx^{-2}e^{-x^2/(2s^2+\sigma^2)}\gg\Pp(Y_1>x)^2$, and hence
\[
 V(x)\ge c x^{-2}e^{-x^2/(2s^2+\sigma^2)}.
\]

For $J_R:=\int_{-R}^R|B_\sigma(x)|^p\rho(x)^{1-p}\,dx$, this variance lower bound yields $\E[|B_\sigma(x)|^p\rho(x)^{1-p}]\ge cx^{-p}e^{\zeta_p x^2}$, where
\[
 \zeta_p=\frac{(p-2)s^2-\sigma^2}{2(s^2+\sigma^2)(2s^2+\sigma^2)}>0
\]
since $\sigma^2<(p-2)s^2$. Thus $\E J_R\to\infty$ as $R\to\infty$. With $m_k=\E|Z|^k$ for $Z\sim\NN(0,1)$, Cauchy--Schwarz gives $\E J_R^2\le(m_{2p}/m_p^2)(\E J_R)^2$, and the Paley--Zygmund inequality implies $\Pp(J_R>\tfrac12\E J_R)\ge m_p^2/(4m_{2p})>0$. Letting $J_\infty=\lim_{R\to\infty}J_R$, we obtain $\Pp(J_\infty>M)\ge m_p^2/(4m_{2p})$ for every $M<\infty$, so $J_\infty=\infty$ with positive probability.

To upgrade this to probability one, define on $\HH_{P,0}$ the integral operator $(\mathcal I h)(x)=\int_{-\infty}^x h(y)\,dy$ and the set
\[
 \mathcal B_p=\left\{h\in\HH_{P,0}:
       \int_\R|\mathcal I h(x)|^p\rho(x)^{1-p}\,dx<\infty\right\}.
\]
By Minkowski's inequality, $\mathcal B_p$ is a linear subspace. By Cauchy--Schwarz, $\|\mathcal I h\|_\infty\le\|h\|_{\HH_P}$, so for each integer $k\ge 1$ the truncated functional $\mathcal J_k(h)=\int_{-k}^k|\mathcal I h(x)|^p\rho(x)^{1-p}\,dx$ is continuous on $\HH_{P,0}$. Consequently,
\[
 \mathcal B_p=\bigcup_{M=1}^\infty\bigcap_{k=1}^\infty
             \{h\in\HH_{P,0}:\mathcal J_k(h)\le M\}
\]
is Borel. By the Gaussian zero--one law for measurable linear subspaces \cite[Corollary~4.56]{HairerSPDE}, $\Pp(G\in\mathcal B_p)\in\{0,1\}$. Since $B_\sigma=\mathcal I G$ and $J_\infty=\infty$ has positive probability, $\Pp(G\in\mathcal B_p)=0$, i.e.\ $J_\infty=\infty$ almost surely.

Finally, on each compact quantile interval $[F(-R),F(R)]$, the inverse-CDF delta method gives
\[
 n^{p/2}\int_{F(-R)}^{F(R)}|F_n^{-1}(u)-F^{-1}(u)|^p\,du
 \Longrightarrow J_R.
\]
Hence $\liminf_{n\to\infty}\Pp(n^{p/2}W_p(P_n,P)^p>M)\ge\Pp(J_R>M)\to 1$ as $R\to\infty$ for every $M<\infty$, proving \eqref{eq:failure}.
\end{proof}

\section{Estimates for unbounded path laws}\label{sec:noncompact-estimates}

This section establishes the weighted transport bounds and feature-moment estimates used in the proofs of Theorems~\ref{thm:scalar} and~\ref{thm:reference}.

\subsection{Quadratic weighted estimates}
\label{sec:optimizedtails}

Take $p=2$ ($\QQ_P=Q_{P,2}$) and enumerate the $N=dT$ scalar coordinates in chronological order $(x_1^{(1)},\dots,x_T^{(d)})$. Write $\bar\rho_j(x_{1:j})=\int\rho(x)\,dx_{j+1:N}$ and $\bar h_j(x_{1:j})=\int h(x)\,dx_{j+1:N}$, with $\bar\rho_0=1$ and $\bar h_0=0$ whenever $\int h=0$. Fix $\beta>0$, $M_0<\infty$, and $B=\beta\sigma^2>\theta_N$ (with $\theta_N$ from \eqref{eq:tailthresholds}), and let $\CC_{M_0}$ denote the class of smoothed densities $\rho=d(\mu*\gamma_\sigma)/dx$ satisfying $\int e^{\beta|a|^2}\,\mu(da)\le M_0$.

\paragraph{Gaussian envelopes and conditional tails.}\label{app:gaussian-envelopes}

Put $\Delta=\alpha-b$ and choose exponents $\alpha,b,\eta,\omega>0$ satisfying
\begin{equation}
 \begin{gathered}
 \alpha>\frac1{2\sigma^2},\qquad
 0<b<\frac{\beta}{1+2\beta\sigma^2},\qquad \eta>0,\\
 \alpha+\eta+2(N-1)\Delta
 <\omega<\frac{\beta}{1+\beta\sigma^2},\qquad \omega<2b.
 \end{gathered}
 \label{eq:exponents}
\end{equation}
The condition $B>\theta_N$ is equivalent to $B^2-(N-\tfrac12)B-(N-\tfrac12)>0$, or
\[
 \frac1{2\sigma^2}+\frac{N-1}{\sigma^2(1+2B)}<\frac{\beta}{1+\beta\sigma^2},
\]
so the feasible range in \eqref{eq:exponents} is nonempty. Set $\HH=\HH_\omega=L^2(\R^N,e^{\omega|x|^2}\,dx)$ and $\HH_0=\{h\in\HH:\int h=0\}$. Integrating $\gamma_\sigma$ over a ball carrying half the mass of $\mu$, completing the square against $e^{\beta|a|^2}$, and using $\E[e^{|x-A_{1:j}|^2/(4\sigma^2)}\mid X_{1:j}=x]\le C e^{\alpha|x|^2}$ yield, uniformly over $\rho\in\CC_{M_0}$ and $1\le j\le N$,
\begin{equation}
 \begin{gathered}
 c e^{-\alpha|x|^2}\le\bar \rho_j(x)\le C e^{-b|x|^2},\\
 |\nabla\log\bar \rho_j(x)|\le C(1+|x|),
 \end{gathered}
 \label{eq:envelopes}
\end{equation}
together with $|\nabla\bar\rho_j(x)|\le Ce^{-b_+|x|^2}$ for any $b_+\in(b,\beta/(1+2B))$.

Write $u=x_{1:j-1}$ and $z=x_j$, with one-step conditional density $\bar k_j(u,z)=\bar\rho_j(u,z)/\bar\rho_{j-1}(u)$ and conditional CDF $\bar K_j(u,z)=\int_{-\infty}^z\bar k_j(u,r)\,dr$. Let $G_c(z)=\int_{|z|}^\infty e^{-cr^2}\,dr$ ($c>0$), set $I_z=(-\infty,z]$ for $z\le 0$ and $I_z=[z,\infty)$ for $z\ge 0$, and put $\Lambda_j(u,z)=\int_{I_z}\bar\rho_j(u,r)\,dr$. Then \eqref{eq:envelopes} gives
\[
 \min(\bar K_j,1-\bar K_j)
 \le \frac{\Lambda_j}{\bar \rho_{j-1}}
 \le C e^{\Delta|u|^2}G_b(z).
\]
For $i<j$, setting $M_{ji}:=\partial_{u_i}\Lambda_j -\Lambda_j\partial_{u_i}\log\bar \rho_{j-1}$, the prefix normalizer $\bar\rho_{j-1}$ cancels from $\partial_{u_i}(\Lambda_j/\bar\rho_{j-1})/\bar k_j$, leaving
\[
 C_{ji}:=\frac{\partial_{u_i}\bar K_j}{\bar k_j}
 =\pm\frac{M_{ji}}{\bar \rho_j}.
\]
Combining this identity with the $b_+$ derivative bound and the score estimate in \eqref{eq:envelopes} yields
\[
 |M_{ji}(u,z)|
 \le C e^{-b|u|^2}G_b(z)
 \le C e^{-b|x_{1:j}|^2},
 \qquad
 |C_{ji}(x_{1:j})|\le C e^{\Delta|x_{1:j}|^2}.
\]

\paragraph{A triangular right inverse and its weighted kinetic energy.}

For $h\in\HH_0$, set $\mathcal I_j(u,z)=\int_{-\infty}^z\bar h_j(u,r)\,dr$ and $B_j=\mathcal I_j-\bar K_j\bar h_{j-1}$, and define $v=(v_1,\dots,v_N)$ recursively by
\begin{equation}
 v_1=-\frac{\mathcal I_1}{\bar \rho_1},\qquad
 v_j=-\frac{B_j}{\bar \rho_j}-\sum_{i<j}v_i C_{ji}.
 \label{eq:scalarrecursion}
\end{equation}
Put $E_j(c)=\int\bar\rho_j(x_{1:j})|v_j(x_{1:j})|^2 e^{c|x_{1:j}|^2}\,dx_{1:j}$ for $c\ge 0$. We show that $\sum_{j=1}^N E_j(\eta)\le C\|h\|_\HH^2$.

On $I_z$, the conditional tail bound above gives $|B_j(u,z)|\le J_j(u,z)+Ce^{\Delta|u|^2}G_b(z)|\bar h_{j-1}(u)|$, where $J_j(u,z):=\int_{I_z}|\bar h_j(u,s)|\,ds$ satisfies $J_j(u,z)^2\le G_\omega(z)\int_{I_z}e^{\omega s^2}|\bar h_j(u,s)|^2\,ds$ by Cauchy--Schwarz. Since $\int_0^\infty e^{(\alpha+c)z^2}G_\omega(z)\,dz<\infty$ for $\alpha+c<\omega$, the direct term $D_j(c):=\int(B_j^2/\bar\rho_j)e^{c|x_{1:j}|^2}\,dx_{1:j}$ obeys $D_j(c)\le C\|h\|_\HH^2$ whenever
\[
 \alpha+c+2\Delta<\omega,\qquad \alpha+c<2b
\]
(with only $\alpha+c<\omega$ needed at $j=1$, since $\bar h_0=0$). For the cross terms ($i<j$), using $C_{ji}=\pm M_{ji}/\bar\rho_j$ eliminates the extra factor $\bar\rho_{j-1}^{-1}$:
\begin{align*}
 \int e^{c|x_{1:j}|^2}\bar\rho_j|v_iC_{ji}|^2\,dx_{1:j}
 &=\int e^{c|x_{1:j}|^2}|v_i|^2
        \frac{|M_{ji}|^2}{\bar\rho_j}\,dx_{1:j}\\
 &\le C\int |v_i(x_{1:i})|^2
             e^{-(2b-\alpha-c)|x_{1:j}|^2}\,dx_{1:j}.
\end{align*}
When $c<2b-\alpha$, integrating out the remaining coordinates $x_{i+1:j}$ and inserting $\bar\rho_i\ge ce^{-\alpha|x_{1:i}|^2}$ yields
\[
 \int e^{c|x_{1:j}|^2}\bar \rho_j|v_iC_{ji}|^2\,dx_{1:j}
 \le C E_i(c+2\Delta).
\]
Now choose $c_j=\eta+2\Delta(N-j)$. For every $i<j$ we have $c_i\ge c_j+2\Delta$, so each cross term at level $j$ is controlled by $E_i(c_i)$, while the direct-term conditions at weights $c_j$ reduce to the strict inequalities in \eqref{eq:exponents}. Induction on $j$ therefore gives $\sum_{j=1}^N E_j(\eta)\le C\|h\|_\HH^2$. Dropping the weight or restricting to $\{|x_{1:j}|>K\}$ yields
\begin{equation}
 \begin{aligned}
 \sum_{j=1}^N\int\bar \rho_j|v_j|^2&\le C\|h\|_\HH^2,\\
 \int_{|x_{1:j}|>K}\bar \rho_j|v_j|^2
 &\le C e^{-\eta K^2}\|h\|_\HH^2.
 \end{aligned}
 \label{eq:energybound}
\end{equation}
For $h\in C_c^\infty\cap\HH_0$, writing $F_i=\bar\rho_{j-1}v_i$ gives $\bar\rho_j v_j=-\mathcal I_j+\bar K_j\bar h_{j-1}-\sum_{i<j}F_i\partial_i\bar K_j$; differentiating in $x_j$ via $\partial_j\bar K_j=\bar k_j$ and the inductive identity $\sum_{i<j}\partial_i F_i=-\bar h_{j-1}$ yields $-\partial_j(\bar\rho_j v_j)-\sum_{i<j}\partial_i(\bar\rho_j v_i)=\bar h_j$, which extends to all $h\in\HH_0$ by \eqref{eq:energybound}.

\paragraph{Metric comparison and continuity of the tangent norm.}

\begin{proposition}[Uniform metric comparison]\label{prop:metric}
For all $\rho,q\in\CC_{M_0}$ and $\omega$ satisfying \eqref{eq:exponents},
\[
 \AW_2(\rho\,dx,q\,dx)\le C\|\rho-q\|_\HH.
\]
\end{proposition}
\begin{proof}
Consider the linear path $\rho_s=(1-s)\rho+sq\in\CC_{M_0}$ ($s\in[0,1]$) with $h=q-\rho\in\HH_0$, and let $R_s:(0,1)^N\to\R^N$ be its inverse Rosenblatt map. For $U\sim\operatorname{Unif}((0,1)^N)$, the coupling $(R_0(U),R_1(U))$ is bicausal and satisfies $\partial_s R_{s,j}(U)=v_{s,h,j}(R_s(U)_{1:j})$. By \eqref{eq:energybound},
\[
 \AW_2(\rho\,dx,q\,dx)^2
 \le\E|R_1(U)-R_0(U)|^2
 \le\int_0^1\sum_{j=1}^N\int\bar\rho_{s,j}|v_{s,h,j}|^2\,ds
 \le C\|h\|_\HH^2. \qedhere
\]
\end{proof}

By \eqref{eq:energybound}, we also have $\QQ_P(h)^{1/2}\le C\|h\|_\HH$ for all $h\in\HH_0$. Setting $\ell=\frac12+\frac{N-1}{1+2B}$, $u=\frac{B}{1+B}$, and $g=u-\ell>0$, an explicit choice of exponents satisfying \eqref{eq:exponents} is $\omega=(\ell+u)/(2\sigma^2)$, $\alpha=(1/2+g/(16N))/\sigma^2$, $b=(B/(1+2B)-g/(16N))/\sigma^2$, and $\eta=g/(8\sigma^2)$.

\subsection{Weighted estimates above the quadratic exponent}
\label{gpstat:sec-high-p}

Now take $p>2$, keeping $\beta>0$, $M_0<\infty$, $B=\beta\sigma^2$, and $\CC_{M_0}$ from Section~\ref{sec:optimizedtails}, and dropping overbars on $\rho_j,h_j,k_j,K_j$ ($\rho_0=1$). Assume $B>\Theta_{p,N}$, where $\Theta_{p,N}$ is defined in \eqref{eq:tailthresholds}.

\paragraph{Parameters, density envelopes, and feature moments.}

Under $B>\Theta_{p,N}$, there exist exponents $\alpha,b,a$ satisfying
\begin{equation}
 \begin{gathered}
 \alpha>\frac1{2\sigma^2},\qquad
 0<b<\frac{\beta}{1+2B},\qquad
 0<a<\frac{\beta}{2(1+B)},\qquad a<b,\\
 p\{\alpha-a+(N-1)(\alpha-b)\}<b.
 \end{gathered}
 \label{gpstat:eq-parameters}
\end{equation}
Indeed, at the boundary values of $\alpha,b,a$, the last inequality becomes $2B^2+\{2-p(N+1)\}B-pN>0$, whose positive root is $\Theta_{p,N}$; evaluating this quadratic at $B=p(N+1)/2$ yields $p>0$, confirming $\Theta_{p,N}<p(N+1)/2$ in \eqref{eq:tailthresholds}.

Fix an integer $s>N/2$ and recall the Hilbert space from Corollary~\ref{gpgeo:hilbertcomparison}:
\begin{equation}
 \HH_{a,s}=\{h:e^{a|x|^2}h\in H^s(\R^N)\},\qquad
 \|h\|_{\HH_{a,s}}=\|e^{a|x|^2}h\|_{H^s},\qquad
 \HH_{a,s,0}=\left\{h\in\HH_{a,s}:\int h=0\right\}.
 \label{gpstat:eq-high-H}
\end{equation}
Sobolev embedding ($s>N/2$) gives $|h(x)|\le C\|h\|_{\HH_{a,s}}e^{-a|x|^2}$, and integrating out $x_{j+1:N}$ yields $|h_j(x_{1:j})|\le C\|h\|_{\HH_{a,s}}e^{-a|x_{1:j}|^2}$, while every $\rho\in\CC_{M_0}$ obeys \eqref{eq:envelopes}. Differentiating $e^{a|x|^2}e_y(x)$ up to order $s$ and completing the square yields, for some integer $L<\infty$,
\[
 \|e_y\|_{\HH_{a,s}}^2
 \le C(1+|y|)^L
 \exp\!\left(\frac{2a}{1-2a\sigma^2}|y|^2\right).
\]
Since $2a/(1-2a\sigma^2)<\beta$ by \eqref{gpstat:eq-parameters}, we obtain
\begin{equation}
 \E_\mu\|e_A\|_{\HH_{a,s}}^2<\infty.
 \label{gpstat:eq-high-second-moment}
\end{equation}

\paragraph{A scalar right inverse with uniform energy tails.}

With $u=x_{1:j-1}$, $z=x_j$, $\Delta=\alpha-b>0$, and $\Lambda_j,M_{ji},C_{ji}$ as in Section~\ref{app:gaussian-envelopes}, we again have $|M_{ji}(u,z)|\le C e^{-b|u|^2}G_b(z)$ and $|C_{ji}(x_{1:j})|\le C e^{\Delta|x_{1:j}|^2}$. For $h\in\HH_{a,s,0}$, set $h_0=0$, $\mathcal I_j(u,z)=\int_{-\infty}^z h_j(u,y)\,dy$, and $B_j=\mathcal I_j-K_j h_{j-1}$. The marginal bound on $I_z$ gives
\[
 |B_j(u,z)|\le C\|h\|_{\HH_{a,s}}
 \bigl\{e^{-a|u|^2}G_a(z)
       +e^{-(a-\Delta)|u|^2}G_b(z)\bigr\}
 \le C\|h\|_{\HH_{a,s}}e^{-(a-\Delta)|x_{1:j}|^2}
\]
(for $j\ge 2$, and without the second term at $j=1$). The recursive velocity $v=(v_1,\dots,v_N)$ defined by $v_1=-\mathcal I_1/\rho_1$ and $v_j=-B_j/\rho_j-\sum_{i<j}v_i C_{ji}$ then satisfies, by induction on $1\le j\le N$,
\[
 |v_j(x_{1:j})|\le C\|h\|_{\HH_{a,s}}
 e^{[\alpha-a+(j-1)\Delta]|x_{1:j}|^2},
\]
since each new coordinate adds at most $\Delta$ to the exponent of the preceding velocities. Consequently,
\[
 \rho_j(x_{1:j})|v_j(x_{1:j})|^p
 \le C\|h\|_{\HH_{a,s}}^p
 e^{-\{b-p[\alpha-a+(j-1)\Delta]\}|x_{1:j}|^2}.
\]
By the strict inequality $p\{\alpha-a+(N-1)\Delta\}<b$ in \eqref{gpstat:eq-parameters}, every stage retains a positive Gaussian decay rate. Choosing $\eta>0$ below the smallest of these rates gives
\begin{equation}
\begin{aligned}
 \sum_{j=1}^N\int |v_j|^p\rho_j
 &\le C\|h\|_{\HH_{a,s}}^p,\\
 \int_{|x_{1:j}|>R}|v_j|^p\rho_j
 &\le C e^{-\eta R^2}\|h\|_{\HH_{a,s}}^p,
\end{aligned}
\label{gpstat:eq-energy}
\end{equation}
and density of $C_c^\infty\cap\HH_{a,s,0}$ in $\HH_{a,s,0}$ gives $h_j=-\sum_{i\le j}\partial_i(\rho_j v_i)$.

\paragraph{Metric comparison and differentiation of quantile maps.}

\begin{proposition}[Uniform weighted comparison]\label{gpstat:prop-comparison}
For all $\rho,q\in\CC_{M_0}$ under \eqref{gpstat:eq-parameters}, the $N$-stage scalar adapted distance $\AW_{p,\mathrm{sc}}$ and the $T$-block distance $\AW_p$ satisfy
\[
 \AW_{p,\mathrm{sc}}(\rho\,dx,q\,dx)
 \le C_{M_0}\|q-\rho\|_{\HH_{a,s}},
 \qquad
 \AW_p(\rho\,dx,q\,dx)
 \le d^{(1/2-1/p)_+}C_{M_0}\|q-\rho\|_{\HH_{a,s}},
\]
and $N_{P,p}$ is continuous on $\HH_{a,s,0}$ for $\rho\in\CC_{M_0}$.
\end{proposition}
\begin{proof}
Along $\rho_\tau=(1-\tau)\rho+\tau q\in\CC_{M_0}$ ($h=q-\rho$, $\tau\in[0,1]$), integrating $\partial_\tau R_{\tau,j}(u)=v_{\tau,h,j}(R_\tau(u)_{1:j})$ on $(0,1)^N$ via Minkowski's inequality and \eqref{gpstat:eq-energy} yields
\[
 \|R_1-R_0\|_{L^p(du;\ell^p)}
 \le\int_0^1\|v_{\tau,h}\circ R_\tau\|_{L^p(du;\ell^p)}\,d\tau
 \le C_{M_0}\|h\|_{\HH_{a,s}},
\]
which proves the scalar bound. The block bound follows from $|z|_2\le d^{(1/2-1/p)_+}|z|_p$ on $\R^d$.
\end{proof}

When $d=1$ and $h\in\DD$, the pointwise velocity bound above and \eqref{gpstat:eq-energy} verify \eqref{eq:scalar-uniform-integrability} with $r=p$. The Vitali argument in Section~\ref{gpstat:sec-scalar-low} then gives $\varepsilon^{-1}(R_\varepsilon-R_0)\to v_{\rho,h}\circ R_0$ in $L^p((0,1)^T;\ell^p)$, completing Proposition~\ref{prop:twopoint} for $p>2$.

\subsection{Strongly log-concave reference laws}\label{sec:referencebounds}

\paragraph{Conditional transport--entropy comparison.}\label{sec:entropy}

Write $\mathrm{KL}(Q\mid P)=\int\log(dQ/dP)\,dQ$ if $Q\ll P$ and $+\infty$ otherwise.

\begin{lemma}[Adapted transport from conditional transport--entropy bounds]
\label{lem:adaptedT2}
Suppose $P$ has finite second moment and positive conditional densities $k_t^P(x_{1:t-1})$ such that for all $1\le t\le T$, $x_{1:t-1},y_{1:t-1}\in\R^{d(t-1)}$, and $\nu\in\mathcal P(\R^d)$,
\begin{align*}
 W_2(\nu,k_t^P(x_{1:t-1}))^2
 &\le 2C_{\mathrm{LS}}\,\mathrm{KL}(\nu\mid k_t^P(x_{1:t-1})),\\
 W_2(k_t^P(x_{1:t-1}),k_t^P(y_{1:t-1}))
 &\le L|x_{1:t-1}-y_{1:t-1}|.
\end{align*}
Then every path law $Q$ with finite second moment satisfies
\[
 \AW_2(P,Q)^2\le C_{\mathrm A}\,\mathrm{KL}(Q\mid P),
 \qquad C_{\mathrm A}=4C_{\mathrm{LS}}(1+2L^2)^{T-1}.
\]
In particular, if $P$ and $Q$ have densities $\rho>0$ and $q$, then
\begin{equation}
 \AW_2(P,Q)^2
 \le C_{\mathrm A}\int\frac{(q-\rho)^2}{\rho}.
 \label{eq:anchoredchisquare}
\end{equation}
\end{lemma}
\begin{proof}
For $\mathrm{KL}(Q\mid P)<\infty$, couple $(X,Y)$ period by period via measurable optimal (or $\varepsilon$-optimal) one-step $W_2$ couplings of $k_t^P(X_{1:t-1})$ and $k_t^Q(Y_{1:t-1})$. Setting $S_t=\sum_{i\le t}\E|X_i-Y_i|^2$ and $H_t=\int Q_{1:t-1}(dy)\,\mathrm{KL}(k_t^Q(y)\mid k_t^P(y))$, the $W_2$ triangle inequality through $k_t^P(Y_{1:t-1})$ gives
\[
 S_t\le(1+2L^2)S_{t-1}+4C_{\mathrm{LS}}H_t,\qquad S_0=0.
\]
Since $\sum_{t=1}^T H_t=\mathrm{KL}(Q\mid P)$ by the chain rule for relative entropy, iterating from $t=1$ to $T$ gives the entropy bound; using $\log(q/\rho)\le q/\rho-1$ then yields \eqref{eq:anchoredchisquare}.
\end{proof}

\paragraph{Conditional constants for the reference class.}

Under Assumption~\ref{ass:structured}, recall $m=\kappa/(1+\kappa\sigma^2)$ and $M=\sigma^{-2}$ from \eqref{eq:Hessianbounds}, and set
\[
 R_f=\frac{C_f}{c_f},\qquad
 D_0=\frac{C_f-c_f}{2\sigma^2\sqrt{\kappa+\sigma^{-2}}}.
\]
The prefix densities of $P$ factor as $\rho_{1:t}=\rho_{0,1:t}r_t$, where $r_t(x_{1:t})=\E_{\mu_0}[f(A)\mid X_{1:t}=x_{1:t}]\in[c_f,C_f]$ satisfies $\|\nabla r_t\|_\infty\le D_0$. At any history $a=x_{1:t-1}$, the conditionals $k_{0,a}$ and $k_a$ under $P_0$ and $P$ obey $\nabla_z^2(-\log k_{0,a})\succeq mI_d$ by \eqref{eq:Hessianbounds}, so $k_{0,a}$ satisfies the logarithmic Sobolev inequality with constant $1/m$ \cite{Caf,FGP}. Because $dk_a/dk_{0,a}(z)=r_t(a,z)/r_{t-1}(a)$ varies by at most a factor $R_f$, Holley--Stroock perturbation and the Otto--Villani theorem \cite{OV,GL} give the conditional $T_2$ inequality for $k_a$ with $C_{\mathrm{LS}}=R_f/m$.

To bound the history Lipschitz constant $L$, take a unit history vector $e$ and decompose
\[
 \partial_e\log\rho_{1:t}(a,z)=F_e(z)+G_e(z),\qquad
 F_e=\partial_e\log\rho_{0,1:t},\quad G_e=\partial_e\log r_t.
\]
From $mI\preceq\nabla^2(-\log\rho_{0,1:t})\preceq MI$, subtracting $\frac12(M+m)I$ shows that the off-diagonal block of $\nabla^2(-\log\rho_{0,1:t})$ between $a$ and $z$ has operator norm at most $M_\times=(M-m)/2$. Hence $|\nabla_zF_e|\le M_\times$, so the Poincar\'e inequality for $k_{0,a}$ gives
\[
 \Var_{k_a}(F_e)\le R_f\Var_{k_{0,a}}(F_e)\le R_f M_\times^2/m,
 \qquad |G_e|\le D_0/c_f.
\]
Since $\partial_e\log k_a$ is the $k_a$-centered version of $F_e+G_e$, it follows that
\[
 \int|\partial_e\log k_a|^2\,dk_a
 \le S^2:=\frac{2R_f M_\times^2}{m}+2\left(\frac{D_0}{c_f}\right)^2.
\]
Interpolating along $a_s=(1-s)a+sb$ and solving $-\nabla_z\cdot(k_{a_s}\nabla u_s)=\partial_s k_{a_s}$, the Poincar\'e inequality for $k_{a_s}$ yields $\int|\nabla u_s|^2\,dk_{a_s}\le(R_f/m)S^2|a-b|^2$, and the Benamou--Brenier formula gives $W_2(k_a,k_b)\le L|a-b|$ with $L^2=(R_f/m)S^2$. Substituting $m,M_\times,D_0$ in terms of $z_0=\kappa\sigma^2$ yields
\begin{equation}
 L^2=\frac{R_f^2}{2z_0^2}+\frac{R_f(R_f-1)^2}{2z_0},
 \qquad
 C_{\mathrm A}=\frac{4R_f}{m}(1+2L^2)^{T-1}.
 \label{eq:explicitentropyconstant}
\end{equation}
Moreover, for any bounded zero-mean tilt $g$ with $|\varepsilon|\|g\|_\infty\le 1/2$, we have $c_f/2\le f(1+\varepsilon g)\le 3C_f/2$ and $\rho+\varepsilon h_g\ge\rho/2$; applying \eqref{eq:anchoredchisquare} with $P_\varepsilon^{h_g}$ as reference law gives
\begin{equation}
 \AW_2(Q,P_\varepsilon^{h_g})
 \le C_*\|q-(\rho+\varepsilon h_g)\|_{\HH_P}
 \label{eq:regularnetcomparison}
\end{equation}
for every probability density $q=dQ/dx$ at finite $\HH_P$-distance.

\paragraph{The intrinsic feature moment.}

\begin{lemma}[Intrinsic second feature moment]\label{lem:intrinsicmoment}
Under Assumption~\ref{ass:structured}, for every $\sigma>0$,
\begin{equation}
 \E_\mu\|e_A\|_{\HH_P}^2
 \le\frac{C_f}{c_f}
      \left(1+\frac{2}{\kappa\sigma^2}\right)^N<\infty.
 \label{eq:intrinsicfeature}
\end{equation}
\end{lemma}
\begin{proof}
Set $q=\mu_0*\gamma_{\sigma/\sqrt 2}$, so that $\rho_0=\mu_0*\gamma_\sigma=q*\gamma_{\sigma/\sqrt 2}$. Applying \eqref{eq:Hessianbounds} at bandwidth $\sigma/\sqrt 2$ shows that $U=-\log q$ obeys $m_*I\preceq\nabla^2 U\preceq 2\sigma^{-2}I$ with $m_*=\kappa/(1+\kappa\sigma^2/2)$. For $Z\sim\NN(0,\sigma^2 I/2)$, Taylor expansion of $U(x-Z)$ gives
\[
 \frac{\rho_0(x)}{q(x)}
 =\E\big[e^{U(x)-U(x-Z)}\big]
 \ge\E\exp\!\left(-\nabla U(x)\cdot Z-\frac{|Z|^2}{\sigma^2}\right)
 =2^{-N/2}\exp\!\left(\frac{\sigma^2|\nabla U(x)|^2}{8}\right).
\]
Since $|\nabla U(x)|\ge m_*|x-x_*|$ around the unique minimizer $x_*$ of $U$, integrating yields $\int(q(x)/\rho_0(x))\,dx\le 2^{N/2}(8\pi/[\sigma^2 m_*^2])^{N/2}$. Using $\E_{\mu_0}[e_A(x)^2]=(4\pi\sigma^2)^{-N/2}q(x)$ and Tonelli's theorem, we obtain
\[
 \E_{\mu_0}\int\frac{e_A(x)^2}{\rho_0(x)}\,dx
 \le\left(\frac{2}{\sigma^2m_*}\right)^N
 =\left(1+\frac{2}{\kappa\sigma^2}\right)^N.
\]
The density bounds \eqref{eq:referenceclass} give $\mu\le C_f\mu_0$ and $\rho\ge c_f\rho_0$, which implies \eqref{eq:intrinsicfeature}.
\end{proof}

\appendix
\section{Spectral formulas and numerical details}\label{app:examples}

\subsection{Quadratic spectrum and a four-point law}\label{sec:spectrum}\label{app:fourpoint}

For $p=2$, write $\QQ_P=Q_{P,2}$ and $\mathcal A_P=\mathcal A_{P,2}$; orthogonal projection onto $\ker(D_\rho)^\perp\subset\mathcal A_P$ defines a bounded linear minimum-energy operator $\mathcal V_P:\HH_0\to\mathcal A_P$.

\begin{corollary}[Trace identity]\label{cor:spectrum}
Under Theorem~\ref{thm:compact}, \ref{thm:scalar}, or~\ref{thm:reference} at $p=2$, the eigenvalues $\{\lambda_j\}$ in Remark~\ref{rem:quadratic} satisfy
\[
 \E\QQ_P(G)=\E_\mu\QQ_P(e_A-\rho)=\sum_j\lambda_j.
\]
\end{corollary}
\begin{proof}
Since $\mathcal V_P:\HH_0\to\mathcal A_P$ is bounded and linear, $\mathcal V_PG$ is a centered $\mathcal A_P$-valued Gaussian element with trace-class covariance $\mathcal V_PC_G\mathcal V_P^*$. Letting $\{u_j\}$ be an orthonormal eigenbasis of $\mathcal V_PC_G\mathcal V_P^*$ in $\mathcal A_P$, the Karhunen--Lo\`eve expansion $\mathcal V_PG\overset d=\sum_j\sqrt{\lambda_j}Z_ju_j$ ($Z_j\sim\NN(0,1)$ i.i.d.) yields Remark~\ref{rem:quadratic} upon taking squared norms, and taking expectations gives $\sum_j\lambda_j = \E\|\mathcal V_PG\|_{\mathcal A_P}^2 = \E_\mu\|\mathcal V_P(e_A-\rho)\|_{\mathcal A_P}^2$.
\end{proof}

The nonzero eigenvalues $\{\lambda_j\}$ also coincide with those of the integral operator $\mathcal K$ on $L^2(\mu)$ defined by
\[
 (\mathcal K g)(a)=\int K(a,b)g(b)\,\mu(db),\qquad
 K(a,b)=\langle\mathcal V_P(e_a-\rho),
                    \mathcal V_P(e_b-\rho)\rangle_{\mathcal A_P}.
\]
Indeed, defining $S:L^2(\mu)\to\mathcal A_P$ by $Sg=\E_\mu[g(A)\mathcal V_P(e_A-\rho)]$ gives $\mathcal V_PC_G\mathcal V_P^*=SS^*$ and $\mathcal K=S^*S$.

For an explicit nonquadratic limit when $p\neq 2$, consider the four-point law $\mu=\operatorname{Unif}\{-a,a\}\otimes \operatorname{Unif}\{-b,b\}$ on $\R^2$ ($a,b>0$). Write $\varphi_\sigma$ and $\Phi_\sigma(x)=\Phi(x/\sigma)$ for the density and CDF of $\NN(0,\sigma^2)$ ($\Phi=\Phi_1$), and put, for $c>0$ and integers $k\ge 1$,
\[
 f_c(x)=\frac{\varphi_\sigma(x-c)+\varphi_\sigma(x+c)}2,\qquad
 D_c(x)=\frac{\Phi_\sigma(x-c)-\Phi_\sigma(x+c)}2,
\]
\[
 s_c(x)=\tanh(cx/\sigma^2),\qquad
 C_{c,p}=\int\frac{|D_c(x)|^p}{f_c(x)^{p-1}}\,dx,
 \qquad \kappa_{c,k}=\int s_c(x)^k f_c(x)\,dx.
\]
Here $\rho(x,z)=f_a(x)f_b(z)$. For any perturbation $h=\rho[u s_a(x)+v s_b(z)+w s_a(x)s_b(z)]$ ($u,v,w\in\R$), formula \eqref{eq:two-scalar} yields
\[
 v_1(x)=-uD_a(x)/f_a(x),\qquad
 v_2(x,z)=-(v+w s_a(x))D_b(z)/f_b(z),
\]
and hence
\[
 Q_{P,p}(h)=C_{a,p}|u|^p+
       C_{b,p}\int|v+w s_a(x)|^p f_a(x)\,dx.
\]
For i.i.d.\ paths $A^{(i)}=(a\varepsilon_i,b\zeta_i)$ with independent Rademacher signs $\varepsilon_i,\zeta_i\in\{\pm 1\}$ and sample averages
\[
 u_n=\frac1n\sum_{i=1}^n\varepsilon_i,\qquad
 v_n=\frac1n\sum_{i=1}^n\zeta_i,\qquad
 w_n=\frac1n\sum_{i=1}^n\varepsilon_i\zeta_i,
\]
the identity $\varphi_\sigma(x\mp a)=f_a(x)[1\pm s_a(x)]$ gives
\[
 \rho_n(x,z)=f_a(x)f_b(z)
 [1+u_ns_a(x)+v_ns_b(z)+w_ns_a(x)s_b(z)].
\]
Since $\sqrt n(u_n,v_n,w_n)\Longrightarrow(Z_1,Z_2,Z_3)\sim\NN(0,I_3)$, Theorem~\ref{thm:compact} gives
\[
 n^{p/2}\AW_p(P_n,P)^p\Longrightarrow
 C_{a,p}|Z_1|^p+
 C_{b,p}\int|Z_2+Z_3s_a(x)|^p f_a(x)\,dx.
\]
By symmetry, $\int s_a(x)^{2m-1}f_a(x)\,dx=0$ for all $m\ge 1$. At $p=2$ the limit reduces to $C_{a,2}Z_1^2+C_{b,2}Z_2^2+\kappa_{a,2}C_{b,2}Z_3^2$, whereas at $p=4$ it equals
\[
 S_4=C_{a,4}Z_1^4+
 C_{b,4}\left(Z_2^4+6\kappa_{a,2}Z_2^2Z_3^2
                          +\kappa_{a,4}Z_3^4\right),
\]
so $n^2\AW_4(P_n,P)^4\Longrightarrow S_4$ and $\sqrt n\,\AW_4(P_n,P)\Longrightarrow S_4^{1/4}$. Replacing $\widehat\mu_n$ by the product of its empirical marginals would replace $w_n$ by $u_n v_n = o_P(n^{-1/2})$, causing the $Z_3$ interaction term to vanish in the limit.

\subsection{Examples with unbounded support}\label{app:gaussian-checks}

Suppose first that $d=1$ and $\mu=\NN(m_0,\Sigma)$ on $\R^T$ with covariance $\Sigma\succeq 0$ (possibly singular). Writing $\lambda_{\max}(\Sigma)$ for the largest eigenvalue of $\Sigma$, we have $\E_\mu[e^{\beta|A|^2}]<\infty$ for every $\beta < (2\lambda_{\max}(\Sigma))^{-1}$, so Theorem~\ref{thm:scalar} applies whenever
\[
 \lambda_{\max}(\Sigma)<
 \begin{cases}
  \sigma^2/(2\theta_T), & 1<p\le2,\\
  \sigma^2/(2\Theta_{p,T}), & p>2.
 \end{cases}
\]
If instead $\Sigma\succ 0$ is nondegenerate in any dimension $d\ge 1$, then $\NN(m_0,\Sigma)$—and any perturbation $\mu=f\NN(m_0,\Sigma)$ with $0<c_f\le f\le C_f<\infty$—satisfies Assumption~\ref{ass:structured} with $\kappa=\lambda_{\max}(\Sigma)^{-1}$. Theorem~\ref{thm:reference} then holds for every $1<p\le 2$ and every fixed $\sigma>0$, with no constraint on $\lambda_{\max}(\Sigma)/\sigma^2$. If the reference potential grows faster than quadratically at infinity, such as $\mu_0(da)\propto e^{-|a|^4-\kappa|a|^2/2}\,da$, then $\E_\mu[e^{\beta|A|^2}]<\infty$ for all $\beta>0$, and Theorem~\ref{thm:reference} holds for all $1<p<\infty$ and all $\sigma>0$.

Under Assumption~\ref{ass:structured}, Lemma~\ref{lem:intrinsicmoment} bounds $I_\mu:=\E_\mu\|e_A\|_{\HH_P}^2$ by \eqref{eq:intrinsicfeature} (with $I_\mu=\det(I_N+\Sigma/\sigma^2)$ when $\mu=\NN(m_0,\Sigma)$). Combining $\AW_2(P_n,P)^2\le C_{\mathrm A}\|\rho_n-\rho\|_{\HH_P}^2$ from \eqref{eq:anchoredchisquare} with $\E\|\rho_n-\rho\|_{\HH_P}^2=n^{-1}(I_\mu-1)$ gives the finite-sample bound
\[
 \E\AW_2(P_n,P)^2
 \le\frac{C_{\mathrm A}}n
 \left[R_f\left(1+\frac2{\kappa\sigma^2}\right)^{dT}-1\right],
\]
where $R_f=C_f/c_f$ and $C_{\mathrm A}$ is given in \eqref{eq:explicitentropyconstant}.

\subsection{Numerical implementation}\label{app:numerics}

For any weight vector $w$ on the four paths in \eqref{eq:experiment-law}, write $P_w=\sum_{i=1}^4 w_i\NN(a_i,\sigma^2I_2)$ and let $R_w$ be its inverse Rosenblatt map. Given the first coordinate $x$, the posterior odds of the upper first-period branch are proportional to $\exp(2ax/\sigma^2)$ and therefore increase with $x$. Conditional on the lower branch, the second coordinate is a Gaussian mixture with means $-2a,0$, whereas on the upper branch its means are $0,2a$; the latter dominates the former in first-order stochastic order for every choice of nonnegative weights. Hence the second-period conditional quantiles of $P_w$ are increasing in $x$ for every $w$, including empirical and bootstrap weights.

For two weight vectors $w,\widetilde w$, let $r(x,u)$ and $s(y,u)$ denote the second-period conditional quantile functions of $P_w$ and $P_{\widetilde w}$. The second-period transport cost is $\int_0^1|r(x,u)-s(y,u)|^p\,du$. Since $r$ and $s$ are increasing in their first arguments and $|r-s|^p$ is submodular, the sum $|x-y|^p+\int_0^1|r(x,u)-s(y,u)|^p\,du$ is submodular in $(x,y)$. The first-period monotone coupling therefore solves the dynamic programming problem (the two-period instance of \cite[Proposition~3.5]{BKR}), yielding
\begin{equation}
 \AW_p(P_w,P_{\widetilde w})^p
 =\int_{(0,1)^2}\sum_{t=1}^2
 |R_{w,t}(u_{1:t})-R_{\widetilde w,t}(u_{1:t})|^p\,du_1du_2.
\label{eq:numerical-kr}
\end{equation}
Each conditional CDF is a mixture of at most three Gaussian CDFs, which we invert by bracketed Newton iteration to CDF residual $10^{-13}$ or bracket width $10^{-12}$, evaluating \eqref{eq:numerical-kr} by tensor-product Gauss--Legendre quadrature with $48$ nodes per coordinate.

To compute the limiting quantiles from \eqref{eq:two-scalar}, put $\psi_i=e_{a_i}-e_{a_4}$ ($i=1,2,3$). The rescaled empirical weight vector $\sqrt n(\widehat w_{1:3}-w_{1:3})$ has covariance $\Sigma_w=\operatorname{diag}(w_{1:3})-w_{1:3}w_{1:3}^{\mathsf T}$, so $G=\sum_{i=1}^3 Z_i \psi_i$ with $Z\sim\NN(0,\Sigma_w)$. Evaluating the three scalar velocity fields from \eqref{eq:two-scalar} with $192$ Gauss--Legendre nodes per coordinate yields a $3\times 3$ energy matrix at $p=2$ (giving the spectrum in Remark~\ref{rem:quadratic}) and a homogeneous quartic polynomial in $Z$ at $p=4$; reference quantiles are estimated from $500{,}000$ independent Gaussian draws.

For the Bermudan put, write $c_i=a_{i,1}$, $b_i=a_{i,2}$, $k=\log(K/S_0)$, and
\[
 H_i=K\Phi\!\left(\frac{k-b_i}{\sigma}\right)
       -S_0e^{b_i+\sigma^2/2}
        \Phi\!\left(\frac{k-b_i-\sigma^2}{\sigma}\right),
\]
the Gaussian expectation of $(K-S_0e^{X_2})_+$ in component $i$. With $\rho_{w,1}(x)=\sum_{i=1}^4 w_i\varphi_\sigma(x-c_i)$, backward induction gives
\[
 \mathsf V(P_w)=\int_\R\max\left\{
 e^{-0.03\Delta}(K-S_0e^x)_+\rho_{w,1}(x),\,
 e^{-0.06\Delta}\sum_{i=1}^4 w_i\varphi_\sigma(x-c_i)H_i\right\}\,dx,
\]
evaluated by adaptive quadrature with absolute and relative tolerance $10^{-9}$ on $[-a-10\sigma,a+10\sigma]$ (omitting Gaussian tail mass worth less than $10^{-19}$ in price units). The unsmoothed value $\mathsf V(\mu)$ is computed by backward induction on the binomial tree.

Each dataset consists of multinomial path counts drawn with probabilities $w$ (random seed $20261116$), and each bootstrap draw resamples multinomial counts from the empirical weights. Using $B_n=249$ draws per dataset, $\widehat q^*_{n,0.95}$ is taken as the $\lceil 0.95\times 249\rceil=237$th order statistic of the scaled bootstrap distances.

\paragraph{Quadrature checks and coverage.}
To verify quadrature stability, we recomputed the first $20$ dataset distances and their first $10$ bootstrap distances at each sample size ($1100$ pairs in total) with $96$ nodes per coordinate; the largest absolute change was $1.2\times 10^{-6}$, and the largest relative change among nonzero distances was $2.2\times 10^{-5}$. Refining the limit calculation from $192$ to $384$ nodes changed the three quadratic eigenvalues—$0.0122484$, $0.0139312$, and $0.0757441$—by less than $3.6\times 10^{-7}$ in relative terms. A small-tilt finite-difference check also agreed with the tangent formula \eqref{eq:two-scalar}. Empirical coverage rates across the $600$ replications are reported in Table~\ref{tab:experiment-coverage}, where each dataset is evaluated with its own bootstrap radius.

\begin{table}[htbp]
\centering
\small
\begin{tabular}{rcccccc}
\toprule
 &\multicolumn{2}{c}{Distance bound}&\multicolumn{2}{c}{Smoothed value}&\multicolumn{2}{c}{Unsmoothed target}\\
$n$&$p=2$&$p=4$&$p=2$&$p=4$&$p=2$&$p=4$\\
\midrule
50&94.2&94.7&100.0&100.0&100.0&100.0\\
200&96.0&96.3&100.0&100.0&100.0&100.0\\
800&94.7&95.5&100.0&100.0&100.0&100.0\\
3200&95.8&96.2&100.0&100.0&57.3&100.0\\
12800&95.2&95.2&100.0&100.0&0.0&0.3\\
\bottomrule
\end{tabular}
\caption{Empirical coverage percentages across $600$ independent datasets at each $n$, using $249$ bootstrap resamples and nominal level $95\%$. Columns report coverage of $\AW_p(P_n,P)\le \widehat q^*_{n,0.95}/\sqrt n$, $|\mathsf V(P_n)-\mathsf V(P)|\le L\widehat q^*_{n,0.95}/\sqrt n$, and the same radius applied around the unsmoothed value $\mathsf V(\mu)$ (Monte Carlo standard error $\approx 0.9\%$ at $95\%$ coverage).}
\label{tab:experiment-coverage}
\end{table}

\section*{Statement on AI use}

Generative AI tools were used during the preparation of this manuscript to assist with proof exploration, intermediate derivations, language editing, and the development of numerical code and figures. The author has independently checked and verified all mathematical arguments, references, code, and numerical results, and takes full responsibility for the content of the paper.

\bibliographystyle{plainurl}
\bibliography{ref}
\end{document}